\documentclass[a4paper,fleqn]{cas-sc}

\usepackage[authoryear,longnamesfirst]{natbib}
\usepackage{caption}
\usepackage{subcaption}
\usepackage{graphicx}
\usepackage[export]{adjustbox}
\usepackage[most]{tcolorbox}
\usepackage{multirow}
\usepackage{float}

\usepackage{etoolbox}
\usepackage{xcolor}

\usepackage{amssymb}
\usepackage{amsthm}
\usepackage{amsfonts}
\usepackage{mathtools}

\newcommand{\taxi}{t_\text{taxi}}
\newcommand{\lifespan}{t_\text{lifespan}}
\newcommand{\timescale}{t_\text{assignment}}

\usepackage{lineno}

\usepackage{booktabs}

\ExplSyntaxOn
\keys_set:nn { stm / mktitle } { nologo }
\ExplSyntaxOff

\ExplSyntaxOn
\cs_gset:Npn \__first_footerline:
  { \group_begin: \small \sffamily \__short_authors: \group_end: }
\ExplSyntaxOff

\def\tsc#1{\csdef{#1}{\textsc{\lowercase{#1}}\xspace}}
\tsc{WGM}
\tsc{QE}

\begin{document}
\let\WriteBookmarks\relax
\def\floatpagepagefraction{1}
\def\textpagefraction{.001}

\shorttitle{Multimodal air-transportation network augmentation}    

\shortauthors{M. M. Borrero and M. Z. Li}  

\title [mode = title]{Augmenting airline networks using airside-to-airside buses to strengthen system resilience under disruptions}

\author[1]{Micah M. Borrero}[orcid=0009-0005-1745-2972]

\cormark[1]

\ead{mborrero@umich.edu}

\credit{Conceptualization, Methodology, Software, Formal Analysis, Investigation, Visualization, Writing – Original Draft.}

\affiliation[1]{organization={University of Michigan},
            city={Ann Arbor},
            state={MI},
            country={United States}}

\author[1]{Max Z. Li}[orcid=0000-0001-8969-6259]

\ead{maxzli@umich.edu}

\credit{Conceptualization, Supervision, Funding Acquisition, Writing – Review \& Editing}

\cortext[1]{Corresponding author}

\begin{abstract}
Each year, disruptions in the air transportation network strand millions of passengers and cost airlines billions in revenue.
Airline networks prioritize operational and cost efficiency through hub-and-spoke structures that maximize revenue; however, these hubs also act as critical choke points during disruptions.
Previous studies focus on \emph{reactionary measures} in response to air transportation network disruptions, whereas this work proposes a \emph{proactive strategy} to improve resilience by reconfiguring the network's topology. 
Specifically, we consider airside-to-airside bus lines as a low-cost, frequent alternative to short, regional flights, offering service that can circumvent air traffic-related delays.
We develop a network construction model that augments the existing air transportation network with these bus lines.
The augmented networks are analyzed through an agent-based simulation, where increased resilience is measured in terms of decreased average hourly passenger delays under both nominal and disrupted conditions.
Our results demonstrate that converting 10 regional routes from air service to airside-to-airside bus service, for a baseline scenario that is constrained by a \$10 million investment budget, can reduce passenger delays by an average of 8\% on disrupted days and 6\% on nominal days.
Furthermore, through a sensitivity analysis, we show that while augmenting the system using these buses decreases operational costs compared to the historic air-only network, continuously expanding bus parameters (such as range and investment budget) ultimately yields diminishing returns in delay mitigation.
Finally, we discuss real-world precedents, alongside regulatory and political hurdles to implementation.
The proposed framework offers airlines, airports, and regulators a decision-support tool for integrating multimodal strategies into future disruption management policies.
\end{abstract}

\begin{highlights}
\item An approach is proposed for multimodal augmentation in an air transportation system
\item Airside-to-airside bus integration reduce passenger delays and operational costs
\item The framework supports preemptive multimodal strategies to improve system resilience
\end{highlights}


\maketitle

\section{Introduction} \label{sec:intro}
Over the course of the 2020s, multiple disruptive events have led to severe air transportation network disruptions. In 2024, Delta Airlines experienced a network disruption instigated by the CrowdStrike outage lasting over five days, with over 7,000 flights canceled~\citep{delta_air_lines_sec_2024}.
In 2022, Southwest Airlines canceled approximately 16,900 flights, leaving 2 million passengers stranded, as a result of Winter Storm Elliott~\citep{us_department_of_transportation_dot_2023}. 

These disruptions are just a couple of events that highlight the susceptibility of the U.S. air transportation network to severe disruptions, where disruptions are categorized as severe if their impacts persist for more than one day~\citep{li_dynamics_2022}.
Disruptions in the air transportation network are inevitable~\citep{marla_integrated_2017, li_dynamics_2022}.
To mitigate the possibility of these disruptions leading to cascades of delays and cancellations, various protocols such as ground delay programs, miles/minutes-in-trail, rerouting, and many other traffic management initiatives (TMIs) have been implemented. These TMIs have enabled in-the-moment handling of disruptive events~\citep{faa_faa_2024, jacquillat_predictive_2022, weitz_nas_2024, clausen_disruption_2010}.
\subsection{Motivation and Research Gap}
While these TMIs provide significant congestion relief during disruptions, they are inherently limited by their reactionary nature. 
In particular, the fundamental hub-and-spoke structure of the present air transportation network acts as a critical congestive bottleneck.
In parallel, multimodal integration has been of interest in recent years, but how this integration impacts the underlying resilience of the air transportation network has yet to be explored in the modern academic literature. 
This presents a significant gap, as, there is a lack of quantitative frameworks designed to evaluate how multimodal alternatives can proactively augment the existing air transportation network to bypass these bottlenecks.
Thus, this research aims to pivot from a reactionary perspective to a proactive approach by directly examining how integrating multimodal alternatives into the present air transportation network impacts the overall future system resilience.
Specifically, we present direct airside-to-airside bus lines as a drop-in alternative to short-haul flights, such as regional flights. We find that bus lines have the potential to serve the dual purpose of improving both \emph{operational resilience} and the \emph{passenger experience}.

\subsection{Overview of Approach}
To address this gap, this study primarily focuses on regional routes in the air transportation network. These routes can be characterized as those where at least one endpoint is a \textit{regional airport}. As there is no universally agreed upon definition of a regional airport~\citep{pauwels_regional_2024}, this work adopts a hybrid definition of what has been proposed by the Federal Aviation Administration concerning the various airport categories~\citep{faa_airport_2022}. In particular, we define a regional airport as one that is either a non-hub (airports that handle less than 0.05\% but more than 10,000 of the annual U.S. commercial enplanements) or a non-primary airport (airports that have between 2,500 and 10,000 annual U.S. commercial enplanements).
These airports are broadly characterized as those with the lowest throughput and connectivity.
Regional airports in a given air transportation network are often serviced by flights that operate at a select number of times throughout the day. These low-frequency flights are primarily a result of airline consolidation to high-value hubs due to relatively low regional demand~\citep{ryerson_integrating_2013, campbell_twenty-five_2012}. Leveraging low-cost, low-capacity buses can improve network connectivity within the air transportation network such that one short-haul flight can be replaced by multiple bus trips. {Consequently, this increased service frequency would help minimize passenger dwell time, ultimately reducing overall trip duration.}
{This is additionally compounded by the fact that these buses are not subject to the typical airspace and airport capacity constraints that constrict air traffic, potentially allowing for continued movement during disruptions.}

Our model identifies regimes where it is both cost-effective and feasible to supplement, if not replace, regional air routes.
This augmented network is constructed leveraging proven methodologies rooted in Linear Programming.
Specifically, we model the air transportation network as a directed graph such that the nodes represent airports, and weighted edges account for the cost to service a given airport pair.
To efficiently allocate edges within the network, we present an integer optimization problem that aims to minimize the cost of each edge throughout the network. 
Utilizing the produced network, we demonstrate the improved resilience of the underlying air transportation network through a case study, which highlights the operational independence and flexibility of these buses. Specifically, we show how these routes are still able to operate, transporting passengers closer, if not to, their intended destination in addition to how the bus fleet can be flexibly leveraged to augment disrupted routes to ensure continued service. Finally, we discuss the present-day applications and policy implications of the integration of airside-to-airside bus lines to the air transportation network.

\section{Literature Review} \label{sec:literature_review}
Various models have been presented that capture the intricacies concerning how the air transportation network responds and recovers from disruptions. Few have considered solutions from a passenger-centric perspective, and fewer have considered multi-modal solutions. In attempting to understand the foundational principles that can be leveraged to develop a given solution, we consider both the qualitative definitions and their quantitative translations, in conjunction with fundamental modeling analogs.

\subsection{Network Disruptions and Resilience}
\label{ssec:netres}
In the context of air transportation research, \cite{janic_modeling_2005} characterized a disruption as a fundamental reduction in the performance of an air transportation network for an extended period of time. 
Simultaneously, when considering a given network's response and recovery from a disruption, one must consider its resilience. 
In the broadest sense, the resilience of a given network can be defined as the ability to sustain operations when disrupted \citep{weitz_nas_2024, gertsbakh_network_2011}.
\citet{hosseini_review_2016} classified various qualitative and quantitative resilience measures, highlighting the broad range of definitions and interpretations of this seemingly fundamental concept. 
Nonetheless, \citet{hosseini_review_2016} emphasized the importance of resilience broadly understood in the general case. 
Focusing on transportation networks,~\citet{miller-hooks_measuring_2012} presented an indicator of resilience rooted in the time it takes to recover from a disruptive event, which has been used in the field of air transportation~\citep{zhang_assessing_2015, janic_reprint_2015}.

The present air transportation literature commonly defines resilience in the context of factors such as delayed flights, disrupted itineraries, or airport capacities~\citep{li_dynamics_2022, barnhart_modeling_2014}.
Although this is understandable from an operational perspective, it falls short in capturing the wide range of possibilities when considering the human element, or more specifically, the passenger. In light of this, \cite{barnhart_modeling_2014} proposed an extensive methodology by which they disaggregated publicly accessible air transportation passenger demand data to model the individual impacts of disruptions on the passenger. They leveraged the passenger delay calculator developed by \cite{bratu_analysis_2005}, which considers factors such as re-accommodation for disrupted itineraries, to estimate passenger delays across the U.S. National Airspace System. For their model, they define disrupted itineraries needing re-accommodation as those such that the non-stop itinerary was canceled or, for a connecting itinerary, the connection would be missed due to significant flight delays. 

Although \citet{barnhart_modeling_2014} did not explicitly propose their methodology in the context of resilience, we can take inspiration from their work and synthesize it with work from the broader field of transportation research. Specifically, we can consider the ``travel time resilience'' throughout the network. Literature in this space predominantly focuses on transit and road networks, but when considering the multimodal component of our work, stark parallels can be drawn such that these concepts can be applied in a broader context. In particular, \citet{mattsson_vulnerability_2015} conducted an extensive review concerning the varying perspectives by which transportation system vulnerability and resilience have been characterized. They described how the series of works, including \citet{jenelius_importance_2006, jenelius_network_2009, jenelius_user_2010}, suggested that we can infer the value of a given link in a network based upon increased travel time required in the event that it is removed. From the perspective of a passenger in an air transportation context, this can be understood as the amount of time lost to a disruption along a given route and the manner by which the passenger weighs the potential time saved by instead taking alternative routes, or itineraries, to reach their intended destination. 

All things considered, it is understandable that this perspective may be brought into question in the present context, considering the degree of autonomy afforded in road networks or personally operated vehicles. With this in mind, we look to work concerning rail and public transportation networks, which typically operate along fixed routes and schedules. In particular, \cite{rodriguez-nunez_measuring_2014} suggested that, in the most basic sense, travelers will often choose the fastest route to reach their destination. Furthermore, in the event that a link, or in our case a route, is closed, they posit that impacted passengers will elect to take the fastest possible detour. In the context of air transportation, this can be viewed as how consumers generally have the highest willingness to pay for non-stop travel itineraries across all U.S. domestic markets---as identified by \cite{ennen_non-stop_2019}. Ultimately, the consideration of passenger travel time in the context of resilience throughout the present air transportation literature has seen minimal applications. 

Nevertheless, while system-level resilience metrics are well-established, the impact of disruptions on passenger travel time has received limited attention in air transportation resilience research. Building on insights from road and rail networks and the observed value passengers place on non-stop itineraries, the present work focuses on quantifying air transportation resilience from a passenger-centric perspective. In particular, we quantify the potential benefits of airside-to-airside bus service as a drop in alternative to regional air service and assess how it can reduce passenger delays under network disruptions.

\subsection{Air Transportation Network Modeling}
\label{ssec:net_model}
Models of air transportation networks have long been formulated as integer programs---as seen through the reviews conducted by~\cite{campbell_integer_1994} and \cite{farahani_hub_2013}. These formulations frequently take the form of Hub Location Problems (HLPs) due to their ability to capture economies of scale on high-density routes and their alignment with the hub-and-spoke structure that emerged in the air transportation industry following deregulation in the 1970s~\citep{campbell_twenty-five_2012, bryan_hub-and-spoke_1999},  where operating cost minimization and network consolidation were greatly incentivized~\citep{ryerson_integrating_2013, shaw_hub_1993, bania_us_1998}. 
Variants of the HLP differ in how they model key network features. Allocation rules may be single allocation---where each spoke is assigned to exactly one hub---or multiple allocation---which allows flows to route via multiple hubs~\citep{campbell_integer_1994, okelly_hub_1994, jaillet_airline_1996, bryan_hub-and-spoke_1999, contreras_branch_2011}. 
Capacity constraints impose upper bounds on hub throughput, shifting models from uncapacitated p-hub problems to capacitated variants~\citep{contreras_stochastic_2011, contreras_branch_2011, abdinnour-helm_hybrid_1998, topcuoglu_solving_2005}

While these models have yielded important theoretical insights and computational tools,  studies concerning both multimodal and air transportation network design have generally focused on the topology of the network in its entirety~\citep{archetti_optimization_2022, farahani_hub_2013}. Although this directly produces an optimal network, this result may not be realistically or foreseeably attainable. Specifically, there may be a variety of hubs that have a degree of prior investment, routes may have regulatory existence requirements that inhibit their adjustment, and ultimately, new structures proposed by these, often, strictly cost driven models may be suboptimal in the context of passenger experience.

With this in mind, we additionally explore models and problems tangential to air transportation network design. There exist various formulations of the Urban Transportation Network Design Problem (UTNDP) that consider the presence of existing network infrastructure and passenger requirements~\citep{farahani_review_2013, duran-micco_survey_2022}. These models provide a framework for designing solutions that effectively integrate real-world constraints, leveraging existing infrastructure while optimizing key factors such as passenger travel time and fleet operating costs.
Although this space is generally well studied, it is noted that there is a large gap between what is done in practice and the application of these models due to data availability and theoretical simplifications~\citep{duran-micco_survey_2022}. 
Furthermore, research in this area frequently emphasizes heuristic-based approaches due to their computational efficiency and adaptability to diverse problem structures~\citep{chau_systematic_2025}. \cite{bertsimas_data-driven_2021} noted that few works have historically addressed this problem using mathematical optimization due to the issue of scaling. Across the various works at scale in this space, most focus on components of a given network that can be modified---such as buses, which can be easily routed as opposed to fixed rail lines. We synthesize the core structure of the HLP to a more restricted network case, incorporating the concepts and techniques presented throughout these various works concerning the UNTDP to reduce the gap between the theoretical models that can be used and real world applications. 
Unifying these discussed models and historical problems, we pivot from the typical approach to strictly consider the presence of an existing network.

Despite the theoretical advances provided by HLP and UTNDP formulations, existing models often fail to account for both operational constraints and passenger-centric considerations in a fully integrated manner. Specifically, while HLPs excel at capturing cost-optimal hub-and-spoke topologies, they rarely incorporate existing network infrastructure or passenger travel time. Conversely, UTNDP approaches better reflect real-world constraints but are limited in their scalability and typically rely on heuristic methods. Building on the strengths of both frameworks, the present work develops a model for an augmented air transportation network that explicitly incorporates pre-existing routes and hubs. This allows for a thorough analysis of the potential resilience benefits that airside-to-airside bus lines might provide under real-world constraints.

\subsection{Multimodal Integration}
Prior work has additionally examined the use of multimodal alternatives to aid in air transportation disruption recovery.
Specifically, both \citet{zhang_real-time_2008} and \citet{marzuoli_improving_2016} have explored the use of charter buses to reroute regional passengers during a disruptive event.
In either study, they proposed optimization and simulation models that jointly determine passenger reassignment across disrupted itineraries and available charter bus capacity. Both works demonstrated that incorporating ground transportation options can substantially reduce total passenger delay and missed connections, particularly during large-scale irregular operations where purely air-based recovery is infeasible.
Additionally, \citet{marzuoli_improving_2016} framed these improvements in terms of operational resilience. They showed that the inclusion of ground transportation expands the set of feasible recovery options and also improves the air transportation system’s ability to maintain connectivity during disruptive scenarios.

In parallel, such multimodal integration for passenger recovery has been explored in the context of Air-Rail systems. \citet{xu_immuner_2023} propose a combined disruption management framework that explicitly leverages high-speed rail services as both a substitute and a complementary mode to air travel. Their optimization model coordinates aircraft and passenger recovery by generating flexible air, rail, and air-rail itineraries. Their approach simultaneously reserves limited air capacity for higher-priority flights and uses available rail capacity to reroute disrupted passengers.
Similarly, \cite{weiszer_air-rail_2025} examine air-rail multimodal operations through a detailed simulation framework that extends an air transport agent-based model, Mercury~\citep{delgado_mercury_2026}, to include rail network processes, passenger transfers, missed connections, rebooking, and station-airport ground access.
By dynamically reallocating passengers across coupled air–rail networks, both works show that rail's integration significantly reduces total delay and improves recovery feasibility and system robustness during highly disruptive events.

While both the bus-focused and the rail-focused literature highlight the significant degree of improvement with respect to sustained passenger movement toward their intended destination, they are separated by the fact that the latter requires a fixed, complementary ground-based network with established intermodal connectivity.
Notably, as noted by \citet{zhang_real-time_2008}, the lack of existing integrated multimodal infrastructure in the United States inhibits the feasibility of these parallel air-rail recovery strategies. 
Conversely, the charter bus-focused concept introduced by \citet{zhang_real-time_2008} and \citet{marzuoli_improving_2016} provide significantly greater operational flexibility, but remain fundamentally limited by its reactionary framework and lack of structured connectivity associated with integrated multimodal networks. 
We intend to bridge this gap by proposing a framework that supports the sustained integration of airside-to-airside bus service. 
Rather than focusing on passenger reallocation during disruption events, this study aims to serve as a proof of concept, investigating how a sufficiently connected ground-based network can be established to support continued system functionality under degraded operating conditions. Furthermore, by operating buses directly within the airside environment, our proposed framework extends multimodal connectivity beyond conventional landside transfer models, enabling stronger integration with existing air transportation operations. 
Through this proposed framework, we seek to combine the operational flexibility provided by charter buses with the structural continuity and network-level integration characteristics of integrated air–rail systems.

\section{Contributions} \label{sec:contribute}

Ideally, existing models could be directly applied to generate optimal network configurations with mode assignments and passenger allocations.
However, in practice, real-world constraints such as the presence of an existing network may prevent these models from providing actionable solutions.
In parallel, a second challenge arises in evaluating the performance of such systems under disruption.
Specifically, standard methods of analysis concerning resilience typically derive metrics in the context of vehicle movements (e.g. delayed aircraft) as opposed to passenger outcomes (e.g., travel delays or missed connections). Although this approach's merits are not disqualified, it falls short with respect to qualifying the true impacts of a given disruption with respect to the passenger.
Considering these presented challenges, the contributions of our work are as follows:

\begin{enumerate}
    \item We develop a stylized model that aids in the identification of cost-effective regimes by which an existing air transportation network may be augmented by intercity, airside-to-airside bus routes.
    \item We propose a definition of resilience rooted in passenger delays and demonstrate the multimodal network's capability to sustain passenger movements during irregular operations, resulting in an 8\% average delay reduction during disrupted conditions.
\end{enumerate}

\section{Methods} \label{sec:method}

\subsection{Demand Data}
\label{ssec:demand-data}
The estimation of demand is a critical component of a given network construction model, as it ultimately determines the resulting network topology. 
To estimate demand markets across the U.S., we primarily utilize historical data from the Bureau of Transportation Statistics' Airline Origin and Destination Survey (DB1B) dataset~\citep{bureau_of_transportation_statistics_transstats_2025}---which is historically a 10\% sample of airline tickets from reporting carriers in the U.S. with respect to each fiscal quarter and is widely used in air travel demand studies for data collected prior to July 2025~\citep{wang_literature_2021, department_of_transportation_office_of_the_secretary_updates_2023}.\footnote{For data reported starting in July 2025, this sample size is increased to 40\% \citep{bureau_of_transportation_statistics_origin-destination_2026}}
For this study, we focus on a typical day during the third quarter of the Fiscal Year 2023. We select this time frame as it falls outside the immediate impacts of the COVID-19 pandemic while still reflecting contemporary passenger demand.

In addition to DB1B, we incorporate the BTS Air Carrier Statistics (T-100) dataset to enhance the accuracy of our model with respect to specific link flow values~\citep{bureau_of_transportation_statistics_transstats_2025}. Unlike the DB1B dataset, which is sampled, the T-100 dataset provides a complete record for each calendar month, reporting all air passenger enplanements for U.S. domestic and international markets, and covers large certified carriers that hold Certificates of Public Convenience and Necessity. Although the dataset offers comprehensive coverage of passenger flows, it lacks detailed passenger origin-destination information. Specifically, using T-100, we can assume a passenger traveled along a given edge, and thus a route must exist between the two airports, but we are unsure as to the true origin and destination of said passenger, whether or not either node was their initial origin or final destination.

Prior to the application of either of these datasets, there is a degree of filtering required.
With respect to the DB1B dataset, we exclude samples where an operating carrier is not reported or where a fare class is not reported, as these are assumed to be ground segments as defined by BTS~\citep{bureau_of_transportation_statistics_transstats_2025}. 
Additionally, entries where ``Market Miles Flown'' are reported as zero are excluded, as these may denote non-air segments or unused tickets.
Furthermore, we exclude samples with more than two hops, as these are deemed to be off-nominal, composing less than 2.5\% of the data. 
For T-100, we exclude entries without scheduled departures, reported passenger counts, or passenger-configured aircraft. 
We keep only carriers in BTS groups 1, 2, and 3---major, national, and large regional U.S. air carriers~\citep{bureau_of_transportation_statistics_transstats_2025, bureau_of_transportation_statistics_number_2024}.

\subsection{Representative Operational Day}
To translate these aggregated passenger volumes into an actionable time frame for our network optimization, we define a representative operational day. 
While the T-100 and DB1B datasets provide quarterly demand samples and monthly passenger volumes, respectively, they do not explicitly capture day-to-day demand variations. 
To address this limitation, we leverage the BTS Airline On-Time Performance dataset to recover the temporal distribution of flights throughout the quarter~\citep{bureau_of_transportation_statistics_transstats_2025}. Using the observed flight schedules and frequencies, we construct a representative operational day that scales the observed historic demand flow between a given route pair by the average fraction of flights operating during active periods. 
This process characterizes typical network activity while preserving historical operational variability.
Notably, the intent is not to reproduce feasible historical schedules, but rather to accurately represent the network's historic topology while accounting for routes that may not have operated daily.
This holistic view enables us to evaluate regimes where it is feasible to augment the existing air transportation network with airside-airside buses.

To ensure no passenger demand is lost in the construction of this representative day, we must also address differences in flight reporting coverage between the T-100 and On-Time Performance datasets.
Specifically, we divide any unmatched flights into two groups: sustained and infrequent operations. 
We characterize sustained service as that where the T-100 dataset indicates periodic flight schedules throughout a given month (i.e., there is an average of two or more scheduled departures performed).
For routes with these recurring flights, we uniformly distribute the aggregate monthly demand into a representative day using the average total number of scheduled departures performed. 
As for infrequent service, we deem these to be off-nominal occurrences. While these may be classified as `Scheduled Passenger/ Cargo Service,' we assume that the lack of consistent scheduling is insufficient to constitute reliable recurring demand for air travel. Thus, these entries are excluded for the purposes of this analysis.

\subsection{Airport and Route Classification}
\label{ssec:class}

\begin{figure}[pos=h]
    \centering
    \includegraphics[width=0.8\linewidth]{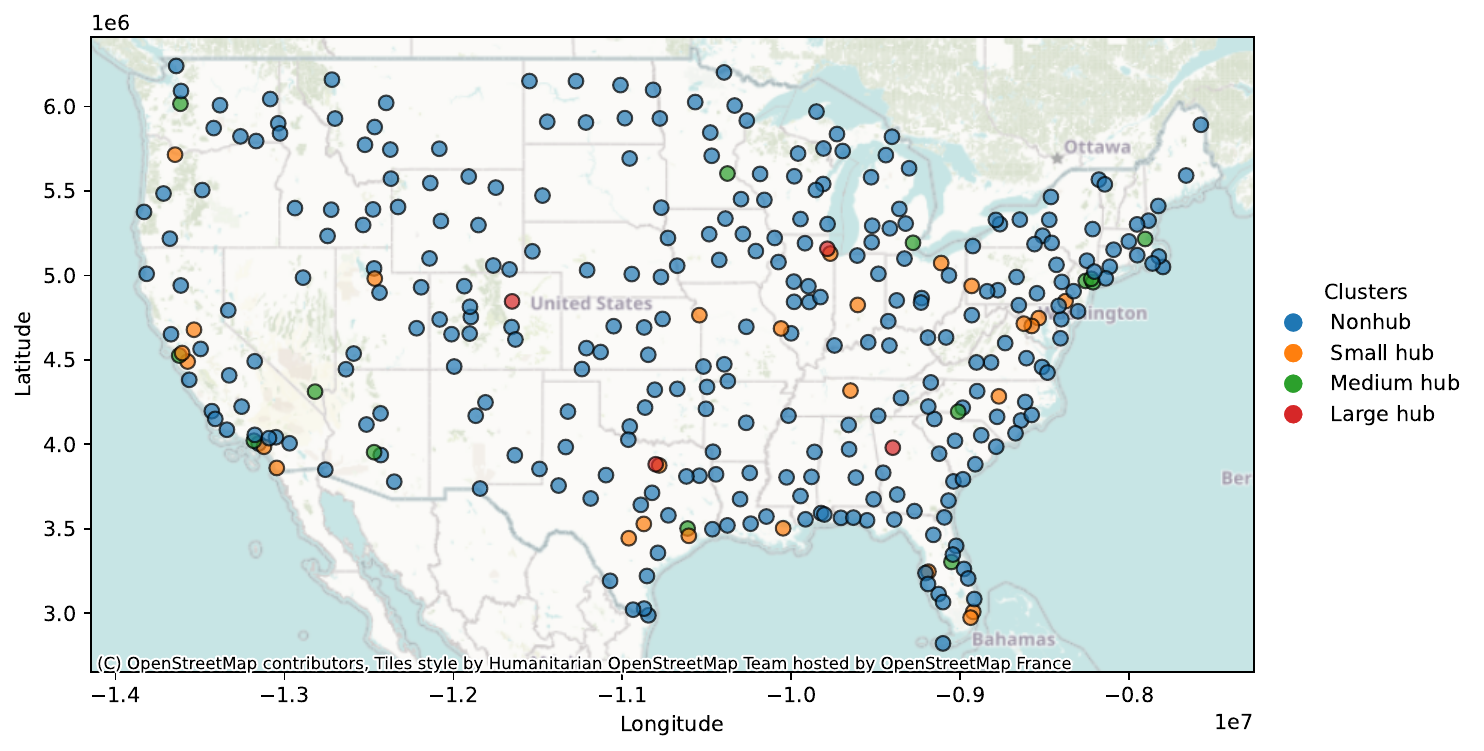}
    \caption{Airport clusters across the continental U.S. for $k=4$.}
    \label{fig:airport-class}
    \vspace{-10pt}
\end{figure}

In this work, we define the historical network graph $G_N = (V_A, E_R)$, where $V_A$ and $E_R$ represent the set of all airports and the directed routes (with respect to vehicle mode) in the network, respectively. 
We assume that, for a given $G_N$, there is a selected portion of the network that must remain constant due to implicit constraints---such as prior investment or modal restrictions.
We define the portion of the network that must remain constant according to a given airport's classification. There are a variety of approaches by which airports may be classified~\citep{pauwels_regional_2024, gao_what_2021, wei_typology_2015}. The primary method used by the Federal Aviation Administration (FAA) classifies airports based on their fractional traffic share over the course of a calendar year or, if it is not a primary airport---a commercial service airport with more than 10,000 annual enplanements~\citep{us_congress_49_2024}---various qualitative and quantitative definitions are used. This may include proximity to a given Metropolitan Statistical Area of a given size for regional airports~\citep{faa_airport_2022}. 
It has been noted that this approach is insufficient in thoroughly characterizing the hub profiles of various airports~
\citep{rodriguez-deniz_classifying_2013}.
As this study aims to capture the spokes within the existing hub-and-spoke air transportation network, the FAA's classification constraints may prove inadequate. For example, the FAA approach may classify an airport with relatively low through-connectivity as a medium or large hub, even though the route servicing it may be classified as a `regional route' by the airline or in the context of node degree. 
With this in mind, we leverage a similar approach to that posited by~\citep{rodriguez-deniz_classifying_2013} to adapt the classification categories defined by the FAA.
Specifically, we classify airports with respect to each node's historic degree and passenger flow, using \textit{k}-means clustering with the squared Euclidean distance metric. 

\begin{figure}[pos=b] %
    \centering
    \includegraphics[width=\linewidth]{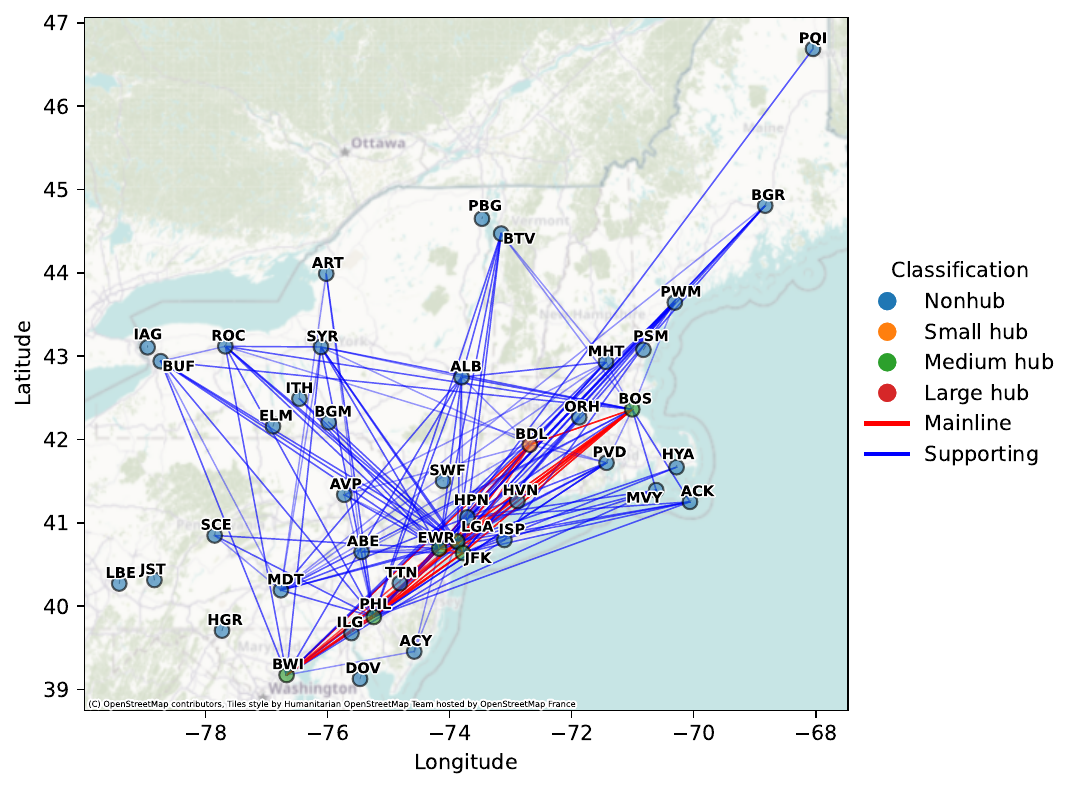}
    \caption{Classified active routes across the Northeastern U.S. Note, as we use $k$-means clustering on historic node degree and passenger flow, traditionally large airports (e.g., PHL, LGA, EWR, JFK) are classified as medium hubs due to their lower domestic through-connectivity compared to airports such as ATL, DEN, and DFW.}
    \label{fig:route-class}
\end{figure}

A node's degree is defined as the number of edges $e \in E_R$
connected to node $v \in V_A$. 
The existing network $G_N$ is constructed using the T-100 dataset
such that a link $e$ exists if it is reported in the dataset by a given air carrier while accounting for the filters described in Sec.~\ref{ssec:demand-data}.
We select $k=4$ in line with commonly used airport classifications, where $\eta_i$ represents the hub classification of airport $i \in V_A$.
These four classes are encoded numerically as $\eta_i \in \{3, 2, 1, 0\}$, corresponding to Large Hub, Medium Hub, Small Hub, and Nonhub, respectively, and are visualized in Fig.~\ref{fig:airport-class}.
We classify routes in $G_N$ as either `mainline' $M$ or `supporting' $S$ routes. 
Mainline routes are defined as those which are fixed, with respect to aviation service, at the time of optimization based upon the values of $\eta_i$ and $\eta_j$ such that the edge connecting city $i$ to $j$ must be serviced by a sufficient number of flights based upon existing demand flows. Specifically, we consider the case where,

\begin{equation}
    \label{eq:mainline}
    M \coloneqq \{(i,j) \in E_R\mid\{\eta_i, \eta_j\} \subset \{2, 3, 4\}\}.
\end{equation}
\\
\noindent From here, it follows that the supporting routes $S$ are defined as the set of all routes in $G_N$ not in $M$.
The set of mainline $M$ and supporting $S$ edges are visualized for the Northeastern U.S. in Fig.~\ref{fig:route-class}.

\subsection{Network Optimization} \label{ssec:optimod}

Before defining the optimization problem, we first define the following. Let $G' = (V_A, E_R')$ be the fully connected graph for all $v \in V_A$ such that $G \subseteq G'$. The following optimization model intends to find the adjusted optimal network graph $G^* = (V_A, E_R^*)$ where $G^* \subseteq G'$. Note that across each subgraph, the vertices under consideration remain constant. This stems from the fact that, given a historical network, all airports that are currently serviced must remain serviced---we cannot elect to exclude a demand set in its entirety.
We leverage the presented classification methodology in conjunction with these definitions and an initial formulation of the hub location problem to define the optimization problem as follows:

\allowdisplaybreaks
\begin{subequations}
    \begin{alignat}{3}
        \min_{y_{i, j}, x_{i, l, j}, a_{i, j}} \quad & \sum_{i, j} \left(\frac {d_{i, j}}{\lambda^{(\alpha)} s^{(\alpha)}} + t_{\text{taxi}} \right)c^{(\alpha)} y_{i, j}^{(\alpha)} + \sum_{i, j} d_{i, j}c^{(\beta)} y_{i, j}^{(\beta)} + \sum_{i, j} \frac {qv_{i, j}^{(\beta)}} {\lifespan} + \sum_{i, j} \frac{ua_{i, j}} {\timescale} \label{eq:obj_func} \\
        \mathrm{s.t.} \qquad~ &\sum_k b^{(k)} y_{i, j}^{(k)} \geq D_{i, j} + \sum_{l \neq i, j} \left( D_{i, l} x_{i, j, l} + D_{l, j}x_{l, i, j} - D_{i, j}x_{i, l, j} \right), && \forall~i \neq j, \label{eq:flow-allocation}\\
        &\sum_l x_{i, l, j} \leq 1, && \forall~l \neq i, j, \label{eq:flow-feasability} \\
        & \left( r^{(k)} - d_{i, j} \right)y_{i, j}^{(k)} \geq 0, && \forall~i \neq j, \label{eq:range-constr} \\
        &b^{(\alpha)} y^{(\alpha)}_{i, j} \geq f_{i, j}^{(M)}, && \forall~(i, j) \in M, \label{eq:mainline-constr} \\
        & \sum_{i, j} qv^{(\beta)}_{i, j} \leq Q, && \forall~(i, j) \in S, \label{eq:bus-aquisition-constr} \\
        &\sum_{i, j} v^{(\alpha)}_{i, j} \leq \sum_{i, j}\bar{v}^{(\alpha)}_{i, j}, && \forall~i \neq j, \label{eq:max-aircraft-constr} \\
        &a_{i, j} \geq v^{(\alpha)}_{i, j} - \bar{v}^{(\alpha)}_{i, j}, && \forall~i \neq j, \label{eq:aircraft-reallocation-constr} \\
        & \frac{1}{2}\left(\frac{d_{i,j}} {\lambda^{(k)} s^{(k)}} + T^{(k)}\right) y_{i,j}^{(k)} \le \mathcal{T} v_{i,j}^{(k)}, && \forall~i\ne j, k, \label{eq:fleet-required-constr}\\
        & v_{i,j}^{(k)} \leq My_{i,j}^{k}, && \forall~i \ne j, k, \label{eq:zero-assign-force}\\
        &y^{(k)}_{i, j} \in \mathbb{Z}_{\geq 0}, &&\forall~i \neq j, k, \label{eq:trip-num} \\
        &x_{i, l, j} \in [0, 1], && \forall~i \neq l \neq j, \label{eq:layover} \\
        &a_{i, j} \in [0, \infty), && \forall~i \neq j, \label{eq:additional}\\
        & v_{i,j}^{(k)} \in \mathbb{Z}_{\ge 0}, && \forall~i\ne j, k. \label{eq:fleet-size}
    \end{alignat}
\end{subequations}

\noindent The model aims to minimize cost $c^{(k)}$ for a given vehicle type $k \in K$ per unit distance $d$, or cost per unit time in the case of the aircraft, where $(i,j,k) \in E_R$. Simultaneously, the model accounts for the lifetime per unit cost of incorporating alternate modes (buses) and the cost of reallocating aircraft throughout the network, scaled according to the estimated assignment length. The objective value produced by this model can be interpreted as the {total operational cost for a representative day.}
We use cost per unit time in the case of the aircraft, as operating costs are often reflected in the form of cost per flight hour. We multiply the cruise speed $s$ by a scalar $\lambda $ to account for reduced speeds during takeoff, landing, and other in-flight maneuvers. 
The addition of $\taxi$ is to account for the average taxi-in and taxi-out times. 

The variable $D_{i,j}$ represents the market level demand between airports $i,j \in V_A$.
$f^{(M)}_{i,j}$ represents the historic passenger flow within the mainline portion of the network $M$. 
The variable $x_{i,l,j}$ represents the fractional flow served by a given intermediary city $l \in V_A$. The set $K$ represents the potential vehicles, or modes, under consideration. We consider the case where a potential fleet is composed of an aircraft and a bus option, defined as the set $K \coloneqq \{\alpha, \beta\}$,
where $\alpha$ represents the aircraft and $\beta$ represents the bus.

\begin{figure}[pos=htb!]
    \centering
    \includegraphics[width=0.6\linewidth]{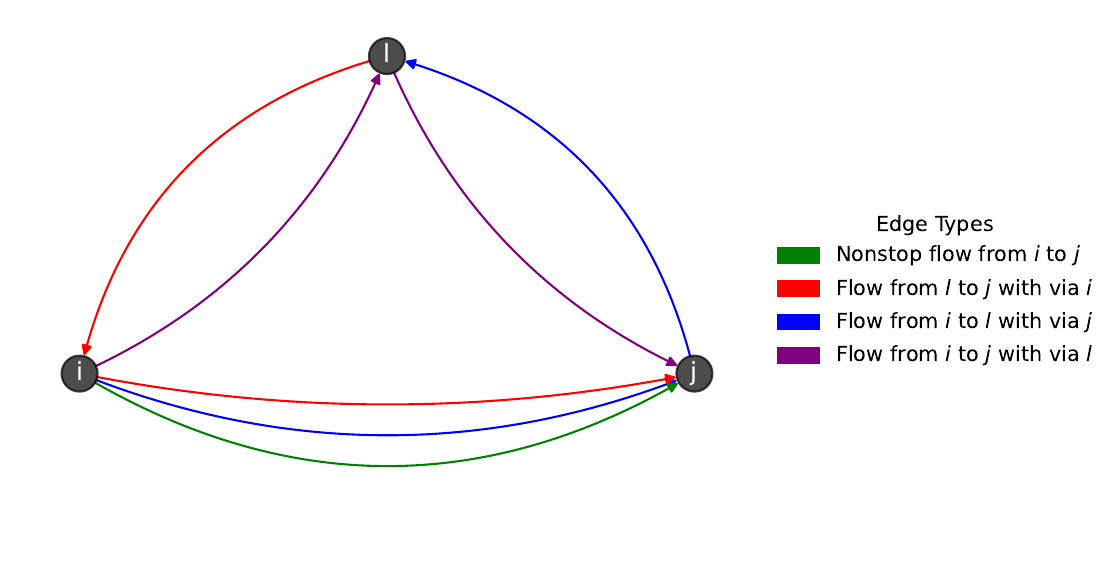}
    \caption{Visualization of passenger flow for a given set of airports $i, l, j \in V_A$.}
    \label{fig:flow-distribution}
\end{figure}

With regard to the presented constraints,~\eqref{eq:flow-allocation} ensures that the total number of trips assigned to a given route is sufficient~\citep{jaillet_airline_1996}. In particular, consider the visualization presented in Fig.~\ref{fig:flow-distribution}. For a given $y_{i,j}^{(k)}$, the number of trips allocated must be able to service the passengers flowing through the red, green, and blue edges during a given operational day. 
This is accounted for in constraint \eqref{eq:fleet-required-constr} where we calculate the number of vehicles $v_{i,j}^{(k)}$ required to successfully complete all trips $y_{i,j}^{(k)}$ for a given route. This is supplemented by constraint~\eqref{eq:zero-assign-force}, which takes the form of a Big-M constraint. Here, $M$ is chosen to be sufficiently large to ensure that if no trips are assigned to a given route pair, then no vehicles are unnecessarily allocated to it.

For the purposes of this analysis, we assume that vehicles are strictly tied to a given O-D pair such that throughout a given operational day, they will only traverse the route between the two nodes. 
While we recognize that, in practice, aircraft and buses may operate with flexible, multi-leg routings, fully integrating fleet scheduling into a multimodal network design problem introduces significant computational complexity. By enforcing this O-D restriction, our optimization model isolates and evaluates the high-level topological feasibility of airside-to-airside bus integration without being confounded by granular fleet routing optimizations.
Furthermore, this approach establishes a conservative baseline for our analysis. We anticipate that, if the augmented network successfully yields resilience benefits and cost efficiencies under rigid, suboptimal fleet assignments, it is highly likely that introducing flexible routing in real-world applications would compound these operational improvements.

With this in mind, we estimate the time to complete a given trip as a function of the vehicle's cruising speed, or highway speed in the case of buses, $s^{(k)}$ scaled by $\lambda^{(k)}$. For aircraft, $\lambda^{(k)}$ accounts for takeoff, landing, and other in-flight maneuvers, whereas for buses it captures speed limitations on access roads. This is combined with a fixed buffer period $T^{(k)}$, which accounts for boarding, disembarkation, and `taxi.' For buses, `taxi' is taken to mean the period spent on airport grounds while not at a stand. Note $\taxi = T^{(\alpha)} \subseteq T^{(k)}$. We divide this by half to account for the potential return trip to the initial point of origin.

With regard to the mainline component, \eqref{eq:mainline-constr} ensures that we are only modifying the supporting component of the network, applying the airport classifications defined in Section~\ref{ssec:class}. 
We assume that for a given historical route serviced by aircraft, such will be operating at, or near, optimal capacity in some form. 
This is incorporated through a frictional cost $u$ applied to each reallocated aircraft $a_{i,j}$,
represented by constraint \eqref{eq:aircraft-reallocation-constr}, where the number of reallocated aircraft is equivalent to the difference in the current number of aircraft assigned to a route minus the historical number,
and accounted for in the third summation within Eq.~\eqref{eq:obj_func}, where we apply the cost of reallocating aircraft resources from one route to another.
This is done in conjunction with the assumption that the size of the aircraft fleet cannot increase due to real-world aircraft acquisition limitations---as represented by constraint~\eqref{eq:max-aircraft-constr}.
We estimate the historic number of aircraft assigned to a given route $\bar{v}^{(\alpha)}_{i,j}$ using a variation of constraint \eqref{eq:fleet-required-constr}. Specifically, we set $k = \alpha$ and $\bar{y}^{(\alpha)}_{i,j}$ is redefined to represent the number of trips historically servicing a given edge estimated by:

\begin{equation}
    \bar{y}^{(\alpha)}_{i,j} = \left\lceil \frac{f^S_{i,j}}{b^{(\alpha)}} \right\rceil,
    \label{eq:historic-trips-max}
\end{equation}

\noindent In this formulation $f^S_{i,j}$ represents the historical flow throughout the supporting portion of the network, $b^{(\alpha)}$ is the capacity of aircraft type $\alpha$, and the term in parentheses in Eq.~\eqref{eq:fleet-required-constr} represents the duration of a single flight segment (inclusive of flight time and a buffer period $T^{\alpha}$). The coefficient $1/2$ is included to account for the return leg to the point of origin, ensuring that the fleet size requirement reflects a balanced round-trip schedule between nodes $i$ and $j$.
We do not directly account for costs or savings associated with a fleet reduction, as this is implied by the first summation in \eqref{eq:obj_func}. Tangentially, we define $q$ as the cost of acquiring a single bus to service a given edge. It follows that we define $Q$ as the total budget earmarked by a fleet operator for the acquisition of buses. 
Constraint \eqref{eq:flow-feasability} ensures that across all intermediary cities, the total fraction of flow serviced by all stopover cities must be feasible. This is visualized as routes such as those represented by the purple edges in Fig.~\ref{fig:flow-distribution}. This prevents vehicle allocations that would otherwise be assigned below the maximal capacity.
Additionally, \eqref{eq:flow-allocation} restricts each trip to a maximum of one stopover. While previous tangential studies have shown that it is possible to improve overall network connectivity and efficiency by increasing the number of stops or transfers---akin to an urban transit system---doing so may decrease convenience and passenger incentivization to travel on these routes, negating any added benefits~\citep{barnhart_modeling_2014}. Thus, we elect to use the one-stop approach.

The constraint defined in \eqref{eq:range-constr} primarily concerns the defined range of a given vehicle such that if the distance between a city pair $(i,j)$ is greater than the serviceable range of a given vehicle, said vehicle will not be allocated to the route directly connecting $i$ to $j$.

\section{ Results and Discussion}
\label{sec:discuss}
To demonstrate the effectiveness of the network generated by the model in Section~\ref{ssec:optimod}, we focus on the Northeastern U.S.---inclusive of New York, New Jersey, and Pennsylvania---as shown in Fig.~\ref{fig:route-class}. In this analysis, we restrict attention to airports with routes entirely contained within the region, namely those connected by the red mainline edges or blue supporting edges in Fig.~\ref{fig:route-class}.

\subsection{Base Case} \label{ssec:base-result}
\begin{table}[pos=h]
    \caption{Baseline model parameters representing typical operational and financial conditions for a Boeing 737 MAX and a standard American motorcoach.}
    \centering
    \begin{tabular*}{\hsize}{@{\extracolsep{\fill}}cll@{}}%
        \toprule
        Variable     & Value                    & Description                            \\ %
        \midrule
        $Q$          & \$10,000,000             & Bus acquisition budget.                 \\ %
        $q$          & \$500,000                & Bus unit cost.                         \\ %
        $u$          & \$10,000                 & Frictional aircraft reallocation cost. \\ %
        $c^{\alpha}$ & \$5,000                  & Aircraft cost per flight hour.         \\ %
        $c^{\beta}$  & \$5                      & Bus cost per mile.                     \\ %
        $b^{\alpha}$ & 160 pax                  & Aircraft passenger capacity.           \\ %
        $b^{\beta}$  & 35 pax                   & Bus passenger capacity.                \\ %
        $r^{\alpha}$ & 3700 mi                  & Aircraft range.                        \\ %
        $r^{\beta}$  & 175 mi                   & Bus range.                             \\ %
        $s^{(\alpha)}$   & 500 mph          & Estimated aircraft max-cruise speed.   \\
        $s^{(\beta)}$    & 55 mph          & Estimated bus highway-speed speed.   \\
        $\lambda^{(\alpha)}$   & 0.8\%          & Conversion factor from aircraft max-cruise speed to average speed.   \\
        $\lambda^{(\beta)}$    & 0.8\%          & Conversion factor from bus highway-speed speed to average speed.   \\
        $\taxi$      & 0.4 hrs                  & Estimated average taxi time. \\
        $\lifespan$  & 7300 days (20 years)     & Estimated vehicle service life.         \\
        $\mathcal{T}$& 17 hours                 & Assumed operating period (0600 to 2300).\\
        $T^{(\beta)}$& 0.25 hours               & Estimated bus `taxi' time.               \\
        \bottomrule
    \end{tabular*}
    \label{tab:nominal}
\end{table}

To ensure that the analysis is consistent, we assume the following nominal values in Tab.~\ref{tab:nominal}.
The values with respect to a given aircraft are estimated based on those of a Boeing 737 MAX due to its relatively broad range of capabilities across short and medium haul routes across the U.S. \citep{boeing_airplane_2025}. Aircraft taxi times were estimated from the ASPM Standard Report of Taxi Times \cite{faa_aspm_2025}.
With regard to the values for the bus, these are estimated based on a standard motorcoach used by some of the largest operators across the U.S., as denoted by the 2020 Motorcoach Census produced by the~\cite{american_bus_association_motorcoach_2025}. In particular, the bus range reflects the maximum distance and inherent travel time that an airline passenger may be willing to spend on a bus. This was determined based upon a combination of the average passenger miles and average passenger trips, scaled to account for a drive time of less than three hours. The bus capacity is estimated based upon the current operating configuration for analogous systems in the Northeastern United States~\citep{perez-castells_even_2024}.

\begin{figure}[pos=htb!]
    \centering
    \includegraphics[width=\linewidth]{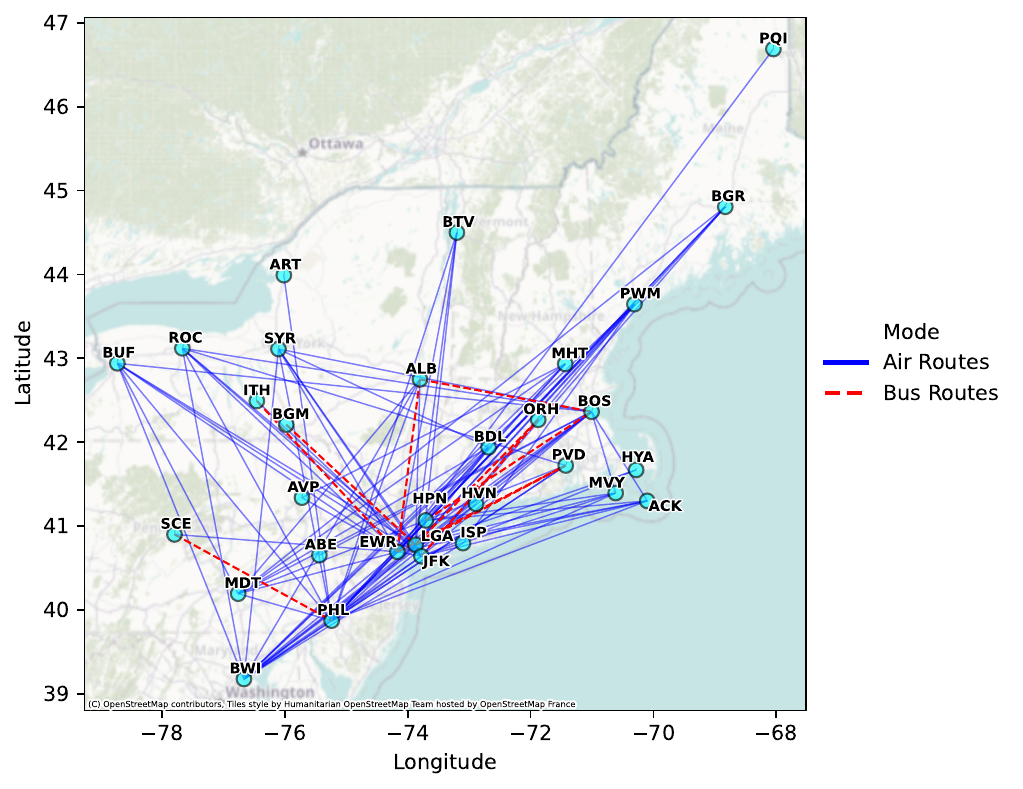}
    \caption{
    Optimization results for an augmented network obtained using the baseline parameter values in Tab.~\ref{tab:nominal}.%
    }
    \label{fig:initial-result}
    \vspace{-40pt}
\end{figure}

To conduct our analysis, we implemented our optimization model in the Gurobi Optimizer version 13.0. For the initial testing of the base case, the solver was run on an Intel Core i7-10 CPU @ 2.60GHz using 6 physical cores and 12 logical processors.
For this initial baseline run, the optimizer converged in under one minute with the results displayed in Fig.~\ref{fig:initial-result}. The resulting objective value, representing the total estimated operational cost for a representative day, was \$8.2 million. For comparison, the estimated operational cost for the historical, air-only network is approximately \$8.34 million. While this is only a 2\% daily cost reduction, the initial \$10 million investment would be recouped within a single fiscal quarter of operations.
Note, although geographical factors are not expressly considered in this analysis, we manually filter bus routes that would be physically infeasible (e.g., routes connecting Nantucket Memorial Airport, ACK, and Martha's Vineyard Airport, MVY).

\subsection{Simulation}
To fully understand the capabilities of a potential air transportation network augmented by bus service, we posit a case study focusing on the Northeastern U.S.---inclusive of New York, New Jersey, and Pennsylvania. A visualization of the routes and airports under consideration can be seen in Fig.~\ref{fig:route-class}.
To simulate the performance of our network, we leverage a simplified agent-based model (ABM) that simulates the flow of aircraft, buses, and passengers throughout the air transportation network. 

ABMs have been used in the Air Transportation Management literature to better understand how system dynamics react to proposed operational improvements and to conduct analysis in controlled environments~\citep{huang_overview_2022, gurtner_agent-based_2021}. Building upon this precedent, we construct our simulation to investigate the potential impacts of an air transportation network augmented by airside-to-airside bus service. Specifically, we simulate the movement of passengers through the selected air transportation network such that delays and disruptions are tracked with respect to the passenger. We subsequently quantify resilience based upon this notion of passenger delay.

\subsubsection{Overview of Scenarios} \label{ssec:sim-scenarios}
To allow for a realistic comparison with respect to the current state of the U.S. National Airspace System, we leverage historical schedule data from the BTS Airline On-Time Performance dataset in this simulation~\citep{bureau_of_transportation_statistics_transstats_2025}.
Using this historic schedule data, we examine four scenarios: three nominal and three disrupted days for both the historic and augmented network.

With regard to the specific days, we utilize the data-driven disruption typologies output by~\citet{xu_identification_2025}. As our study primarily focuses on the Northeastern U.S., we draw from their ``Northeastern Disruption'' cluster. This category is defined by the severity and geographic distribution of delays and cancellations across the NAS, explicitly capturing days with significant operational breakdowns concentrated in this region of interest.
From the outputs within this cluster, we selected a sequence of three heavily disrupted days: 14, 15, and 16 July 2023.
To establish a baseline while controlling for day-of-week effects, and mitigating the broader impacts of temporal or seasonal variability, we pair these disrupted days with three nominal days from an adjacent week: 21, 22, and 23 July 2023.

\subsubsection{Network Construction}
For either base case, the network structure under examination is drawn directly from the historical data as observed across the various TransStats datasets (e.g., DB1B, T-100, etc.), such that if a flight is observed, its corresponding edge is assumed to exist. That being said, we restrict our analysis to the set of airports defined at the beginning of Sec.~\ref{sec:discuss} and visualized in Fig.~\ref{fig:initial-result}. Additionally, we rely on the filters defined in Sec.~\ref{ssec:demand-data} to ensure consistency with respect to the results produced by the presented optimization model. As the network structure predominantly mirrors that observed in the historic schedule data, we assume that the historic schedule sufficiently captures the demand given some latent Average Load Factor (ALF).

With regard to the augmented case, the network structure is drawn from the optimization results discussed in Sec.~\ref{ssec:base-result}. 
As a historic schedule does not exist in its entirety for this augmented network, we adapt the data as follows: For trips scheduled on aircraft edges, we assume that these remain. We assume that the number of historic flights must be greater than or equal to the number of trips estimated by the optimization model. This is due to the fact that the model assumed that each flight must be at or near capacity---whereas in reality, we know that this may not always be the case.
With this in mind, if there are more scheduled flights than the number of trips we reduce the number of scheduled flights while accounting for the estimated Average Load Factor (ALF) and the distribution of the time of day demand.
Based on historical data sourced from the Bureau of Transportation Statistics and the International Air Transportation Association (IATA)~\citep{bureau_of_transportation_statistics_transstats_2025, international_air_transport_association_economics_2025}, we estimate the ALF in recent years to be approximately 83\%---excluding the period between 2020 and 2022 due to the impacts of the Covid-19 pandemic. This estimation is combined with the distribution of peak departure times, which typically occur in the early morning and early evening~\citep{brey_latent_2011}. As such, we proceed with the assumption that flights outside of these prime windows would be the first to be cut in the presence of extreme excess capacity. 

\subsubsection{Schedule Adaptation} \label{ssec:sim-sched}
With respect to bus scheduling, when a bus is introduced to replace a route segment that was historically operated by an aircraft, and the associated tail number was exclusively assigned to flights serving that segment, we perform a direct 1:1 substitution. In this case, the tail number is replaced by a synthetic identifier with a `B-prefix' to denote a bus-operated route, the capacity is fixed to 35 seats (as defined in Tab~\ref{tab:nominal}). In cases where the historical tail serviced multiple routes, we first determine whether these additional routes have been designated as bus edges in the augmented network. If so, we again apply a 1:1 substitution, replacing the aircraft assignment with a bus assignment. If this is not the case, but the historical tail has exactly one remaining leg not assigned to a bus, then the tail is restricted to that single leg, and the other legs are reassigned to a bus with a new B-prefix identifier. Finally, if multiple legs remain, we sequentially evaluate the remaining flights to determine whether the destination of one flight aligns with the origin of the subsequent flight. When this connectivity condition is satisfied, the same substitution procedure described above is applied. In this analysis, no additional cases of this type were observed and, therefore, no further rules were implemented.

If no historical schedule is available for a given route, we construct one using the demand data described in Sec.~\ref{ssec:demand-data}. Beginning with 0800 as the reference time, flights or bus trips are scheduled around this peak morning period. The number of trips is scaled according to demand and the assigned vehicle type: when demand exceeds the capacity of a single vehicle, additional trips are scheduled in adjacent hourly blocks before and after the peak. Service is limited to the operating window of 0600 to 2300. To minimize the fleet size required on each route, we account for the time needed for a vehicle to complete a round trip and return to its origin before assigning the next departure. Consistent with the fleet-sizing approach defined in constraint \eqref{eq:fleet-required-constr} of the optimization model in Sec.~\ref{ssec:optimod}, if the round-trip duration exceeds the operational day, the number of vehicles assigned to the route is increased, and the scheduling process is repeated.

\subsubsection{Capacity Assignments and Scaling}
To extract the passenger capacity of each aircraft, we leverage the FAA's Aircraft Registration Database in conjunction with the OpenSky Aircraft Database to identify the associated aircraft model type for each tail number reported in the filtered Airline On-Time Performance Dataset~\citep{faa_aircraft_2025, schafer_bringing_2014}. For simplicity, we standardize each aircraft's capacity based on its average seat number for a given aircraft type and tune individual values when explicitly known.

When the number of flights simulated within the region of interest is reduced, the representation of capacity-induced delays may become biased. To mitigate this effect, the estimated arrival and departure capacities at each airport are scaled by the proportion of simulated flights relative to the total daily flight count. For instance, if airport ABC has an arrival and departure capacity of 10 flights per hour, but only 6 of the 12 scheduled flights through ABC are included in the simulation, the adjusted capacity is set to 5 flights per hour.

\subsubsection{General Assumptions}
In this simulation, we primarily account for aircraft delays in the context of airport capacity. This allows us to generalize the sources of capacity-reducing disruptions---e.g, adverse weather. Recognizing that a significant portion of observed delays may not be directly related to airport capacity, or the reduced capacity level may be unknown for a specific day, true historic departure delays are enforced to capture the impacts of latent disruptors.
As this simulation's purpose is to aid in the examination of delays as impacted by environmental constraints, we assume that various operational factors, such as crew scheduling or the need for reserve crews, are not a source of delay. 
Additionally, we do not explicitly model schedule recovery efforts to both reduce overall model complexity and obtain a proper baseline for the final result.

\subsubsection{Simulation Specifications}
Independent of the assumptions described previously, we define a set of modeling assumptions that strike a balance between realism and functionality. 
We define three types of agents in this simulation: Airport, Vehicle, Flight, and Passenger.

\paragraph{Airport Agents}
Airport agents define the environment in which the system operates. Each airport is characterized by its 3-letter code, departure and arrival capacities, and queues of passengers waiting for specific flights. Capacities are dynamically reset at the top of each hour to reflect time-varying operational constraints.
Airports manage the arrival and departure of aircraft while also allowing buses to bypass air-based constraints.
Passenger queues are maintained for each scheduled flight, and the airport releases passengers to board flights according to vehicle capacity and boarding rules.

\paragraph{Vehicle Agents}
Vehicle agents represent aircraft and buses operating throughout the network and are identified by their N-Numbers (tail numbers) and B-Numbers (bus numbers) respectively. 
Each agent tracks its current location, status (waiting, turnaround, arrived), scheduled sequence of flights, and time spent in turnaround between flights. 
Vehicle agents automatically advance through their flight schedule and transition between statuses according to arrival and departure rules at airports, in addition to turnaround completion. Vehicles designated as buses are flagged automatically by their B-Number, which alters boarding and capacity rules for associated flights.

\paragraph{Flight Agents} 
Flight agents represent both scheduled aircraft and bus `flights.'
These agents are characterized by their flight number, origin, destination, scheduled departure time, true departure time (if applicable), scheduled flight duration, assigned vehicle, current status (waiting, boarded, enroute, arrived), accrued delay, and elapsed time enroute. 
Boarding occurs no later than 15 minutes before scheduled departure, contingent on vehicle availability and airport capacity. Once a passenger agent is onboard, the flight agent updates both its own and the passengers' status throughout the journey. Departure is conditioned on the vehicle being at the correct location, and the airport allowing departure. If the vehicle is unable to depart once boarded, or land once in the vicinity of its destination, it then accumulates delay. 

\paragraph{Passenger Agents}
Passenger agents represent units of moving demand throughout the network. Each passenger stores its itinerary, current leg, current location, status (waiting, queued, in-flight, missed-connection, arrived), current vehicle, and accrued delay. Passengers queue at airports for their assigned flight and are released when the flight boards. They update their status and delay based on the progression of the flight, subject to the constraints imposed on other agents.

\paragraph{Simulation Execution}
The simulation tracks time in minutes, hours, and days starting from a user designated start time. Each simulation step represents one minute---where all agents execute their individual step functions, updating their statuses, locations, and interactions according to the rules described above. This structure allows for the realistic propagation of delays in the context of the passenger while also considering network congestion effects.

\subsubsection{Simulation Results}
\begin{figure}[pos=htb!]
    \centering
    \includegraphics[width=0.8\linewidth]{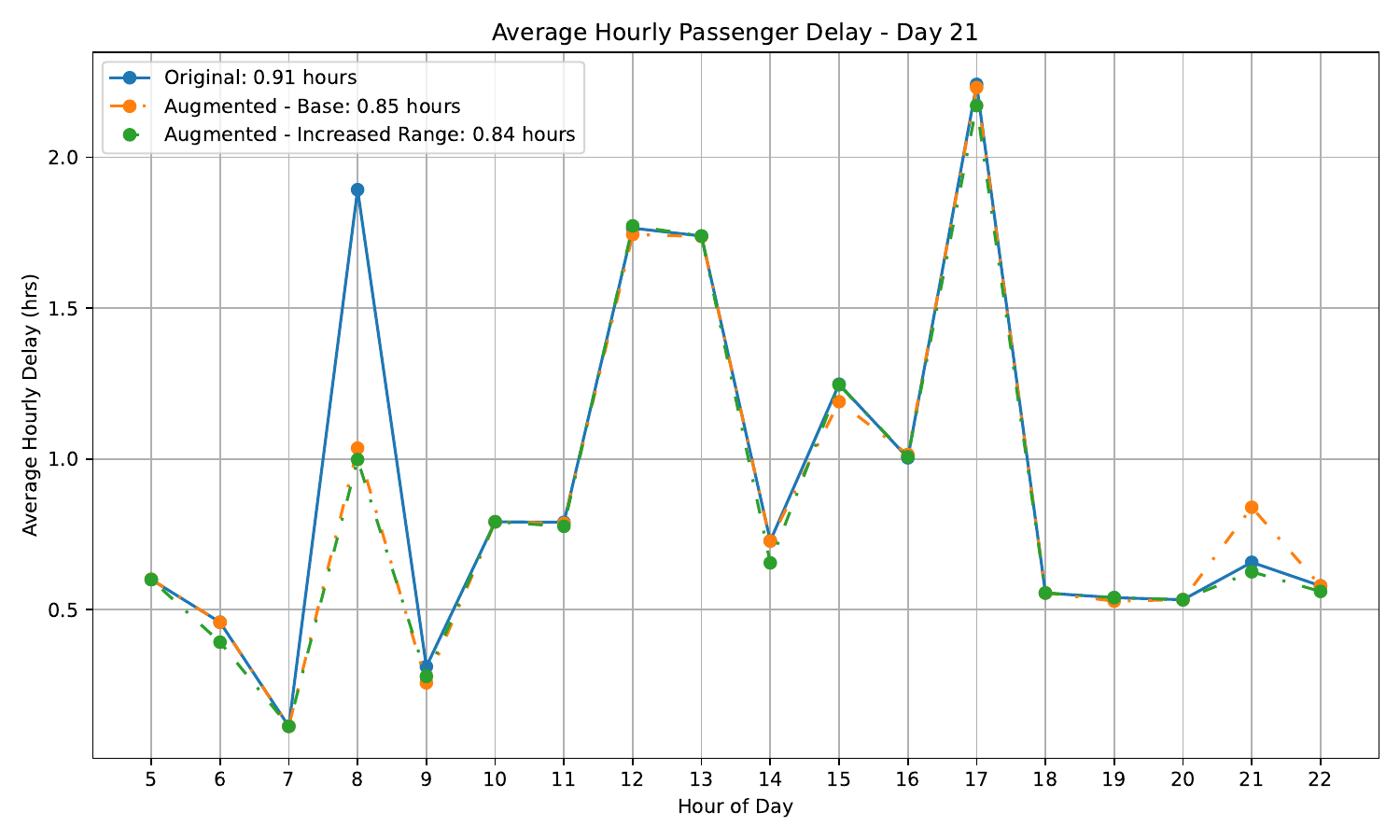}
    \caption{Average hourly passenger delay for both the original and augmented New England network(s) for the nominal case on July 21\textsuperscript{st}, 2023.}
    \label{fig:sim-results-21}
\end{figure}

\begin{figure}[pos=htb!]
    \centering
    \includegraphics[width=0.8\linewidth]{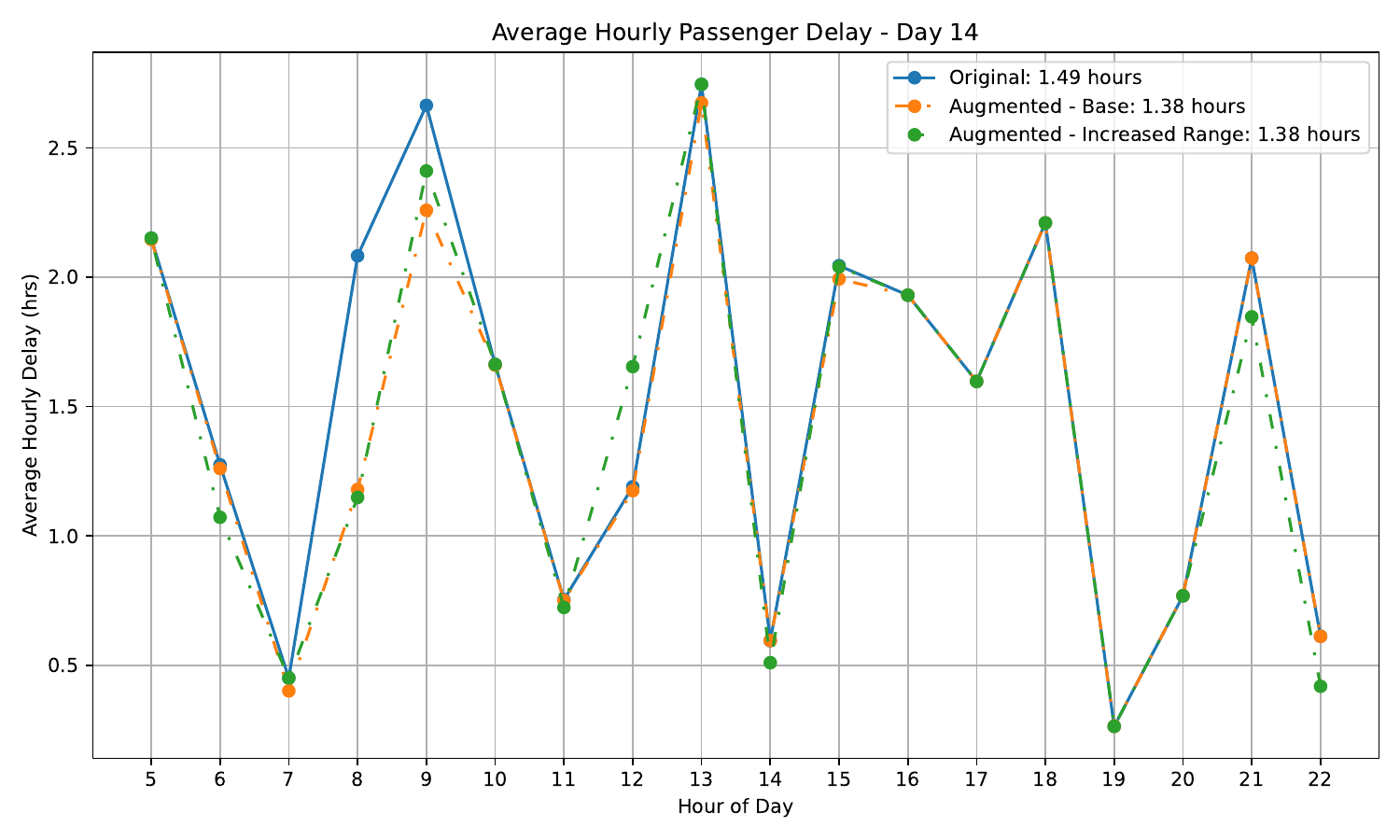}
    \caption{Average hourly passenger delay for both the original and augmented New England network(s) for the disrupted case on July 14\textsuperscript{th}, 2023.}
    \label{fig:sim-results-14}
\end{figure}

To evaluate the four scenarios described above, we ran the ABM simulation across six calendar days (three nominal and three disrupted). Table~\ref{tab:abmsimres} reports the average hourly passenger delay for each date. For a more granular view, Fig.~\ref{fig:sim-results-21} and Fig.~\ref{fig:sim-results-14} illustrate the hour-by-hour breakdown for a nominal day and a disrupted day, respectively. The granular visualizations for the additional four days are provided in \ref{apdx:addl-viz}.

\begin{table}[pos=h]
    \caption{Simulation results showing average passenger delay values, calculated as the mean across hourly observations for each day.}
    \centering
    \begin{tabular*}{\hsize}{@{\extracolsep{\fill}}ccccc@{}}
        \toprule
        Date & Avg Pax Delay - Original (hrs) & Avg Pax Delay - Augmented (hrs) & \% Reduction & Day Classification\\ 
        \midrule
        14 July 2023 & 1.49 & 1.38 & 7.4\% & Disrupted \\ 
        15 July 2023 & 0.75 & 0.70 & 6.7\% & Disrupted \\ 
        16 July 2023 & 1.33 & 1.19 & 10.5\% & Disrupted \\ 
        21 July 2023 & 0.91 & 0.85 & 6.6\% & Nominal \\ 
        22 July 2023 & 0.49 & 0.47 & 4.1\% & Nominal \\ 
        23 July 2023 & 0.58 & 0.55 & 5.2\% & Nominal \\ 
        \bottomrule
    \end{tabular*}
    \label{tab:abmsimres}
\end{table}

In Fig.~\ref{fig:sim-results-21}, the average hourly passenger delay on July 21, 2023, is estimated to be 0.91 hours (55 minutes) for the baseline (original) network. With the introduction of the augmented network, this value is reduced by approximately $6.6\%$, resulting in an estimated average delay of 0.85 hours (51 minutes). Similarly, Fig.~\ref{fig:sim-results-14} shows that on July 14, 2023, the average hourly passenger delay in the baseline network is 1.49 hours (89 minutes). Under the augmented configuration, delays are reduced by approximately $7.4\%$, yielding an average of 1.38 hours (83 minutes).
Note, in both Fig.~\ref{fig:sim-results-21} and Fig.~\ref{fig:sim-results-14} we also display the simulated result for an augmented network for buses with an increased maximum range. This will be explored in detail in Sec.~\ref{ssec:sense}.

Interestingly, in Tab.~\ref{tab:abmsimres} we observe a relatively high average hourly delay across the days classified as `nominal.' These elevated values are a direct consequence of the simulation's design structure rather than an exact reflection of real-world passenger experience. As detailed previously, historical departure delays are strictly enforced to capture latent network disruptors, yet real-world schedule recovery efforts---such as reserve crew deployment, or ad-hoc aircraft swaps---are not explicitly modeled. 
As a result, once a delay is injected into the schedule, it persistently propagates and accumulates across subsequent flight legs without mitigation. 
While this lack of simulated recovery mechanisms inflates the observed absolute delay values, it does so equally across both the original and augmented networks. As a result we appropriately observe higher delay values overall for disrupted days compared to nominal days. 

Additionally, in Fig.~\ref{fig:sim-results-21} we observe a brief period, around 2100, where delays are increased for the augmented network, compared to the baseline network. 
This is most likely attributed to the way the schedule is constructed, where certain flights and aircraft are replaced by buses early in the day, leading to downstream schedule impacts.
Specifically, as the augmented network replaces these aircraft with synthetic bus schedules heavily anchored to the morning peak, the evening period likely suffers a localized capacity shortfall. 
Furthermore, the rigid 35-seat bus capacity means that if evening passenger demand exceeds these limited evening schedules, or if earlier bus trips fall behind, passengers accrue delay as they wait for a subsequent trip. This spillover effect, combined with reduced evening capacity, leads to the accumulation of delays in the later hours.

Overall, we observe that the introduction of airside-to-airside bus services contributes to delay reductions in both the nominal and disrupted cases, with more pronounced benefits observed in the disrupted case as observed in Tab.~\ref{tab:abmsimres}. 
In both Fig.~\ref{fig:sim-results-21} and Fig.~\ref{fig:sim-results-14}, these reductions are concentrated around 0800. Similar to the source of the observed evening delays in the augmented case, this concentration stems from the design of the bus schedule, which was constructed to align with the morning departure bank, as described in Sec.~\ref{ssec:sim-sched}. As a result, the observed delay reductions are temporally clustered around this block. During this period, delay reductions reach approximately \emph{$40\%$} in the nominal case and \emph{$43\%$} in the disrupted case. 

\subsection{Sensitivity Analysis} \label{ssec:sense}
\begin{figure} [pos=htb!]
    \centering
    \includegraphics[width=\linewidth]{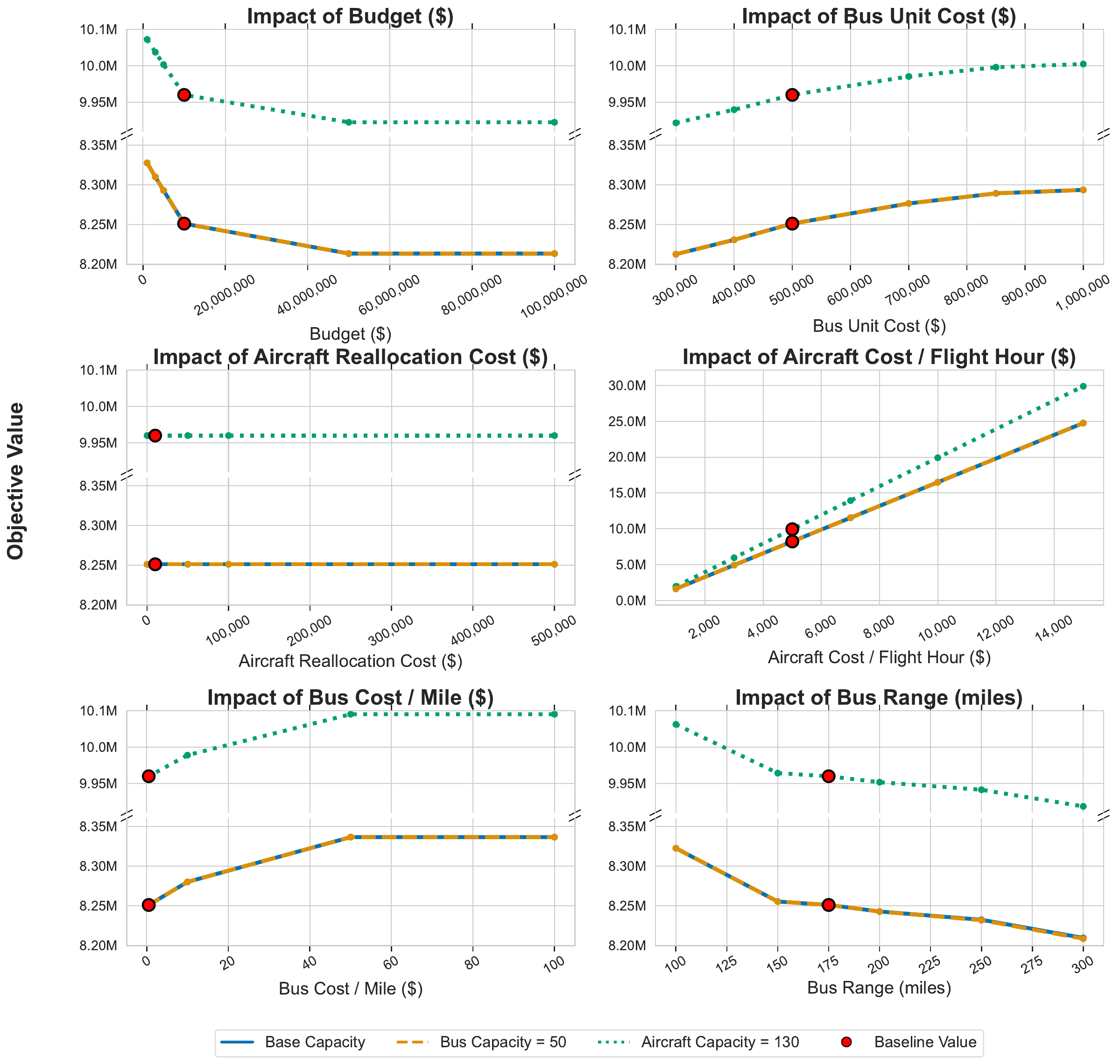}
    \caption{Sensitivity with respect to objective value (total estimated cost for a given operational day).}
    \label{fig:obj-sense}
\end{figure}

The following subsection presents a sensitivity analysis of the model's output to critical input assumptions.
All constant input values are drawn from Tab.~\ref{tab:nominal}.
For a given air carrier exploring a potential investment in bus lines as an alternative to regional air service, it is assumed that funds are not simply unlimited. Furthermore, variables such as operating costs, bus range, fleet reallocation costs, investment budget, and vehicle passenger capacity  may differ from our estimated values.
Thus, to enhance our understanding of the limitations and extent of an ideal investment, we examine the 
impacts of varied input parameters on the resulting objective value, or total daily operational cost, produced by our model. The input values that we vary include investment budget, potential unit cost, reallocation costs, operating cost, and vehicle range. In parallel to these parameters, we examine two alternative vehicle capacities for the bus and aircraft, this is to account for the fact that carriers frequently 
operate diverse vehicle configurations and non-homogeneous fleets based upon localized needs.

\begin{figure} [pos=htb!]
    \centering
    \includegraphics[width=\linewidth]{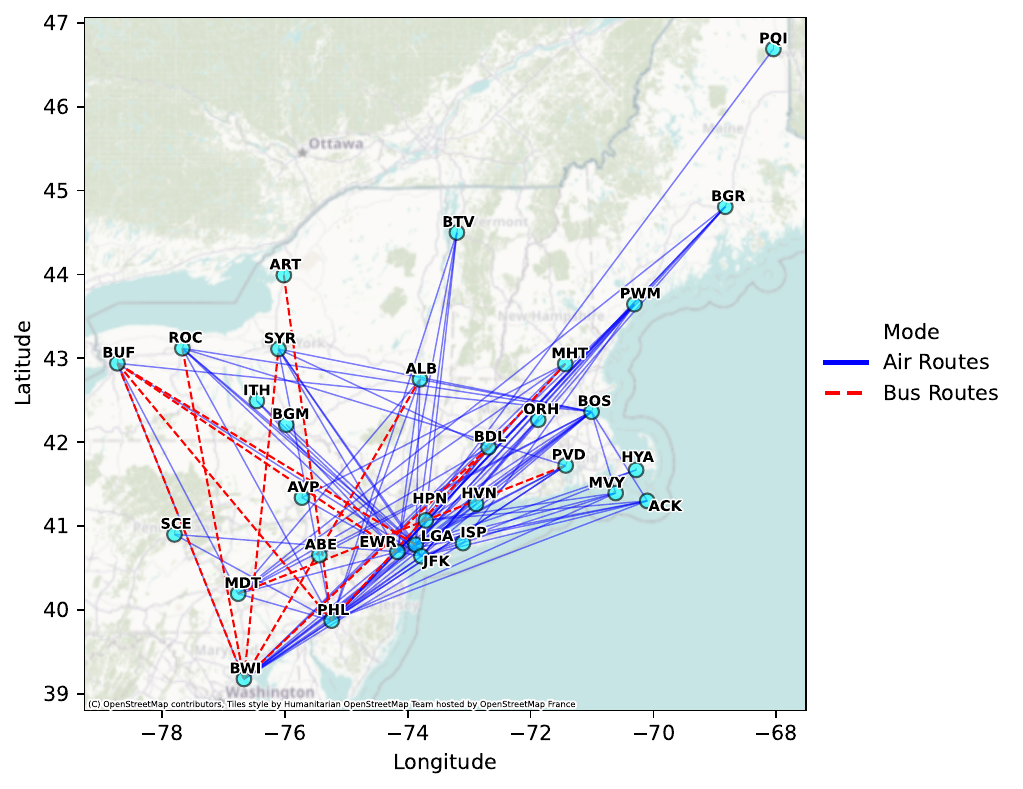}
    \caption{Optimization results for an augmented network considering buses with an increased range of 300 mi.}
    \label{fig:lr-bus}
\end{figure}

In Fig.~\ref{fig:obj-sense}, we present the results of this sensitivity analysis, where, certain linear trends emerge. For instance, we observe a positive relationship between increases in aircraft operating cost per flight hour and the objective value. Additionally, we identify an inverse relationship between the objective value and increases in bus range. 
While the former positive relationship is trivial, the latter trend can be explained by the comparatively lower baseline operating cost per mile of buses relative to aircraft, combined with an extended range leading to the network visualized in Fig.~\ref{fig:lr-bus}. Beyond these intuitive patterns, other parameters, such as bus expenses, investment budget, aircraft reallocation costs, and vehicle capacity, exhibit more nuanced behaviors.

The analysis reveals that aircraft reallocation costs have a negligible impact on the objective value. This can be attributed to the fact that reallocation occurs infrequently and its costs are high compared to the daily aircraft operating expenses. For bus operating costs, the effect plateaus after a threshold. This suggests that, although operating costs rise, they remain relatively lower than those of aircraft, and further increases fail to meaningfully shift the cost tradeoff. Network topology may also play a role, limiting the benefits of additional investment once a certain allocation level has been reached. A similar plateau effect is observed with the investment budget: the objective value decreases sharply at first, then stabilizes. This indicates that, beyond a certain threshold, the maximum beneficial number of bus routes has already been assigned, given the other system parameters.

We additionally observe that varying the bus capacity has minimal impact on the resulting objective value. This is likely due to the fact that, in the current model, these buses are not operating at their maximum passenger limits---indicating potentially lower load factors. Conversely, decreasing aircraft capacity correlates with an increase in operational costs. This behavior is broadly to be expected; on routes where it may not be feasible to assign a bus, serving equivalent demand levels necessitates a larger fleet to accommodate a higher flight frequency. Interestingly, the objective value scales at a higher rate for the lower capacity aircraft with respect to increases in aircraft operating costs per flight hour. This can be attributed to a compounding effect of requiring a larger fleet and higher flight frequencies which amplifies the network's overall sensitivity to fluctuations in aviation operational expenses.

Pivoting to the operational impacts of these parameter variations, we examine the simulated schedule discussed in Sec.~\ref{ssec:sim-sched}. In particular, we simulate the schedule for the case where the maximum bus range is increased to 300 mi, producing the previously mentioned network visualized in Fig.~\ref{fig:lr-bus}. For the disrupted case in Fig.~\ref{fig:sim-results-14}, we observe that this change has virtually no impact on the average passenger delays across the operational day. 
While there are slight hour-to-hour fluctuations relative to the augmented base case, these variations ultimately prove to be  negligible. For the nominal case in Fig.~\ref{fig:sim-results-21}, we do observe a slight improvement compared to the base augmentation case with a reduction in the average passenger delays from 0.85 to 0.84 hours. These two factors combined with the objective value for this scenario indicate that while increasing the serviceable route distance for buses may lower overall operational costs, its impact on passenger delay yields diminishing returns.

Overall, this analysis reveals that augmenting the air transportation network with airside-to-airside bus service lowers total daily operational costs. However, as aircraft operations continue to account for the majority of operating expenses, the magnitude of this cost reduction remains relatively small. 
Furthermore, we find that while increasing various input parameters---such as bus range, or the overall investment budget---can lower estimated operational costs, these modifications may yield diminishing returns with regard to passenger delay mitigation. Although not all parameter variations were discussed in this analysis, further investigation is required to explore the holistic impacts of these parameter variations on schedule performance.
Despite these limitations, we find that this multimodal integration improves the underlying flexibility and connectivity of the network while simultaneously serving as a cost-effective strategy to enhance system resilience.

\section{Current Practice and Policy Needs} \label{ssec:policy}
\begin{figure} [pos=htb!]
    \centering
    \setlength{\fboxsep}{0pt}%
    \fbox{\includegraphics[width=0.75\linewidth]{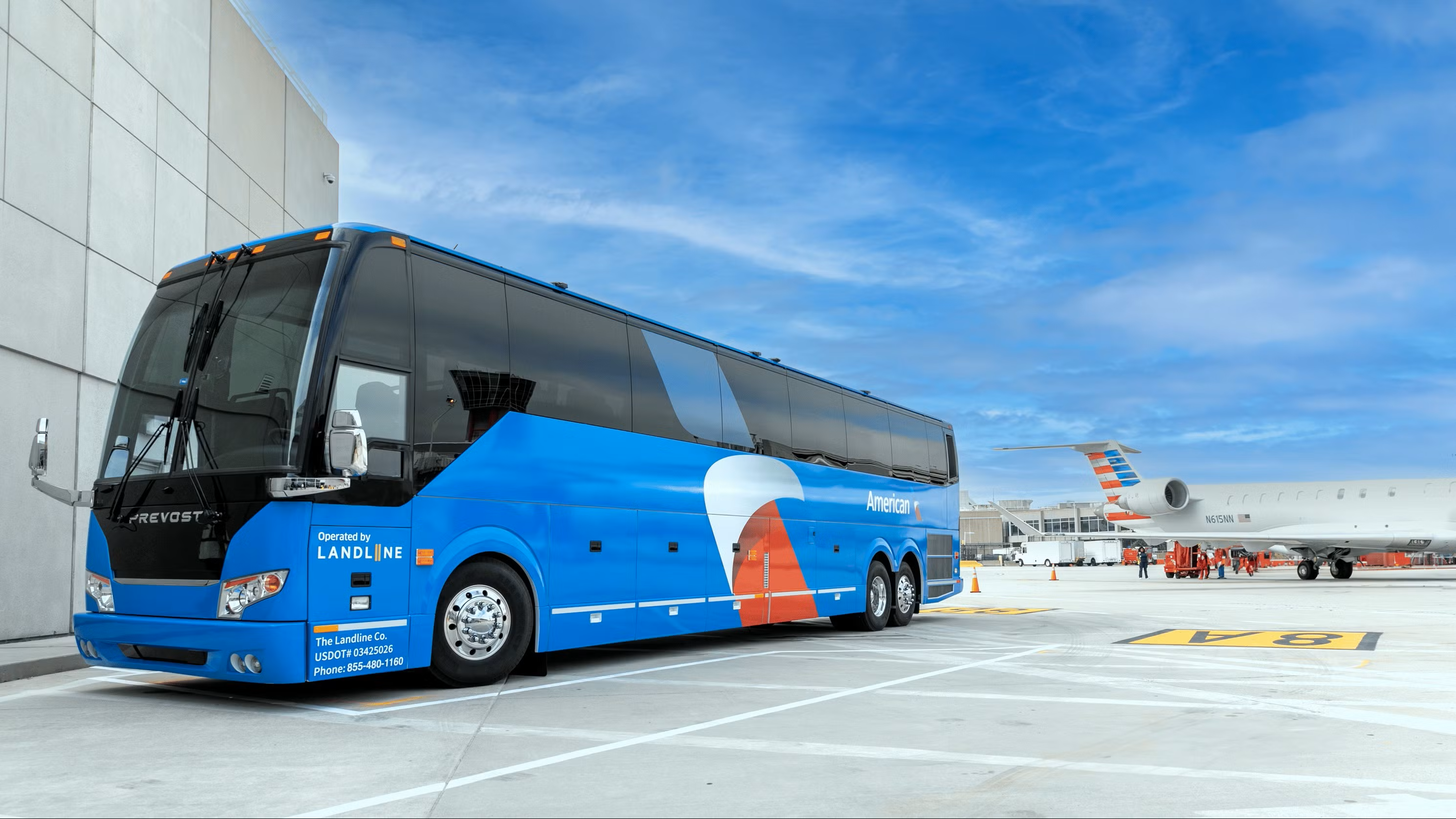}}
    \caption{Landline Company partnership with American Airlines~\citep{hartley_american_2025}.}
    \label{fig:landline-bus}
\end{figure}

The concept described in this work is already seeing real-world applications. In North America, The Landline Company has partnered with various North American air carriers, such as American Airlines, Sun Country, and Air Canada, to integrate both landside-to-landside and airside-to-airside bus service~\citep{landline_company_landline_2025}, with examples of these partnerships shown in Fig.~\ref{fig:landline-bus}. Similar partnerships have been developed in Europe with the Lufthansa Express Bus in Germany and the KLM Bus in the Netherlands \citep{lufthansa_lufthansa_2025, klm_royal_dutch_airlines_travel_2025}. 
These multimodal partnerships demonstrate how this air-bus integration is not simply a conceptual exercise, but a rapidly evolving practice in the air transportation industry. 

\begin{figure*}[pos=htb!]
    \centering
    
    \begin{subfigure}[t]{0.48\textwidth}
        \centering
        \fbox{\includegraphics[width=2.75in,height=2.75in]{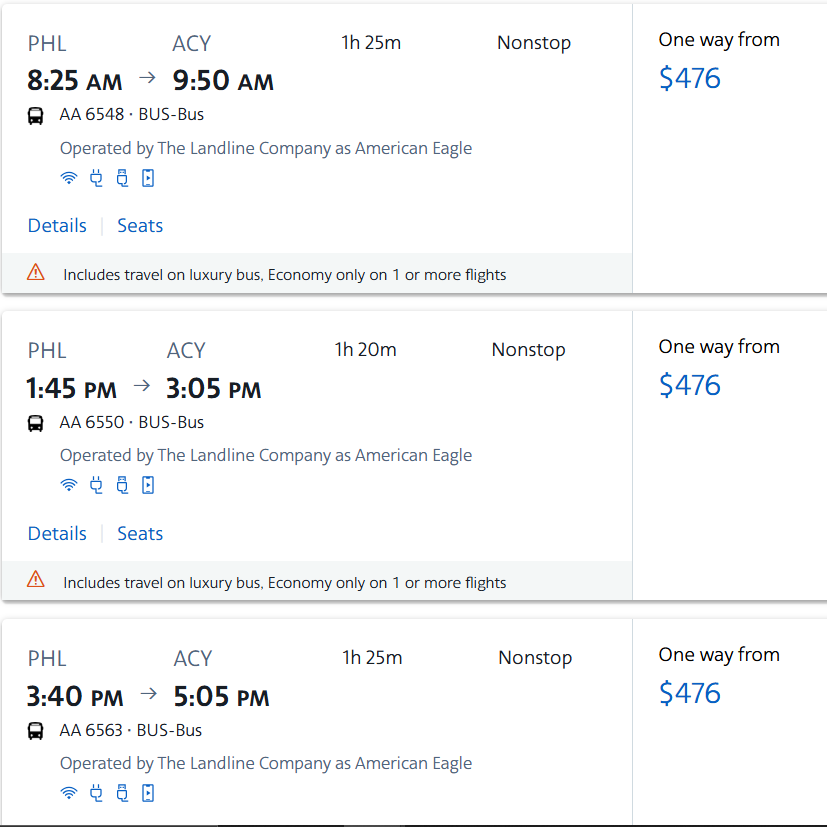}}
        \caption{Screenshot of American Airlines booking interface for a PHL–ACY itinerary on Sept 30, 2025 (captured on Sept 23, 2025)~\citep{american_airlines_american_2025}.}
        \label{sfig:american}
    \end{subfigure}
    ~
    \begin{subfigure}[t]{0.48\textwidth}
        \centering
        \fbox{\includegraphics[width=2.75in,height=2.75in]{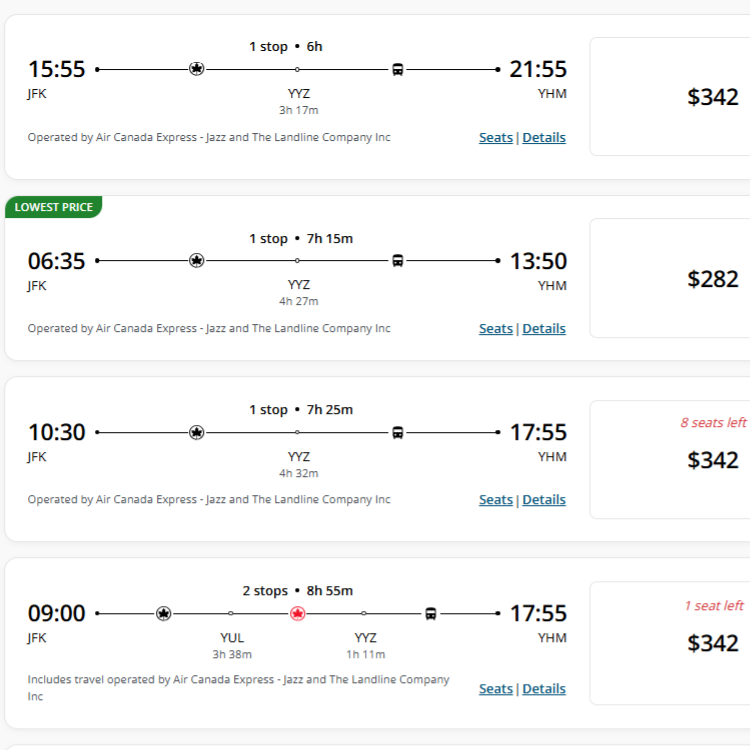}}
        \caption{Screenshot of Air Canada booking interface for a JFK–YHM itinerary on Sept 30, 2025 (captured on Sept 23, 2025)~\citep{air_canada_air_2025}.}
        \label{sfig:canada}
    \end{subfigure}

    \caption{Sample screenshots of airline booking options on buses.}
    \label{fig:ticket-examples}
\end{figure*}

The examples shown in Fig.~\ref{fig:ticket-examples} highlight how seamlessly these multimodal options have been integrated, expanding regional service and connecting passengers to a broader set of regional destinations. 
Nonetheless, a closer examination reveals that the current implementations remain limited in scope. Specifically, the available itineraries are primarily designed to facilitate connections on broader, longer-haul trips rather than serving as standalone regional links. For instance, in Fig.~\ref{sfig:american}, the targeted route for service by American Airlines between Philadelphia (PHL) and Atlantic City (ACY) is effectively priced out as an independent segment, limiting its accessibility for passengers seeking shorter regional connections. Similarly, in the Air Canada example shown in Fig.~\ref{sfig:canada}, the leg from Toronto (YYZ) to Hamilton (YHM) is not offered as a standalone option; it only appears within the context of a longer connecting itinerary from YVR to another major hub. This suggests that while the current multimodal offerings demonstrate technical feasibility, further development is needed to fully leverage bus-to-air integration as a flexible and widely accessible regional mobility solution. 

These limitations highlight the importance of policy frameworks and institutional guidance in shaping the adoption and impact of multimodal service. Airports and airlines, being the primary stakeholders, must navigate regulatory, operational, and financial considerations to expand these services effectively. The following sections examine these policy implications, focusing on how airport infrastructure, airline network strategy, and regulatory oversight influence the feasibility and scalability of bus-air integration.

\subsection{Airports}
From an operational standpoint, present applications have shown the minimal impact that the addition of airside bus service has on airport operations due to its minimal footprint when compared to aircraft. This limited impact makes airside buses an attractive tool for managing passenger flows, particularly at airports where gate and airfield capacity are constrained. At the same time, the substitution of bus segments for short-haul flights introduces financial implications. As airports generate a large portion of their revenue through landing fees, terminal space rentals, and Passenger Facility Charges (PFCs), a reduction in aircraft movements may result in lost revenue when considering the current profit model. Nonetheless, new revenue opportunities may emerge as airports develop dedicated agreements for airside bus operations, positioning themselves as multimodal hubs rather than strict gateways to the sky. 

Balancing these trade-offs will require airports to evaluate whether the strategic benefits of reduced congestion and expanded regional connectivity outweigh the potential revenue loss from fewer aircraft operations. This cost-benefit analysis is particularly critical for smaller regional airports that rely on Federal support through the Essential Air Service (EAS). In particular, present regulatory structures for funding are tied to `enplanements' such that if these are reduced, an airport's eligibility for funding may be reduced~\citep{us_department_of_transportation_essential_2024}. That being said, in 2003, the Federal government established a pilot program known as the Alternate Essential Air Service (AEAS)~\citep{us_department_of_transportation_alternate_2019}. The AEAS was meant to act as a supplement to the core EAS, providing communities with opportunities for funding to support alternative options, such as regularly scheduled ground transportation. To our knowledge, the AEAS has seen limited application beyond the pilot phase. Nevertheless, as the transportation system continues to evolve, programs such as the AEAS can serve as a base to adapt broader programs, such as the EAS, for multimodal service.

\subsection{Airlines}
For airlines, the policy implications of bus-air integration are equally multifaceted. On the one hand, replacing marginal short-haul flights with airside bus service can reduce operating costs, improve network resilience, and expand connectivity into smaller regional markets without the expense of maintaining thinly served routes. This strategy has already been employed in practice, where airlines ticket and brand bus segments as part of their official network, as seen in the partnerships described at the beginning of Sec.~\ref{ssec:policy}. On the other hand, airlines must weigh how the incorporation of contracted bus service will impact competitive positioning, network structure, and regulatory compliance. 

Specifically, airlines must determine how bus operations are incorporated into alliance agreements and minimum service obligations under the EAS. 
Particularly when considering the EAS, regulatory approval may be needed to substitute bus service for air service while still receiving federal support. 
Additionally, although airlines may benefit financially from outsourcing thin routes, they must also manage passenger perceptions. In particular, passengers in North America may often overlook buses due to perceived inconvenience compared to driving or flying. Nevertheless, as noted by \cite{ryerson_drive_2018} and \cite{barnhart_modeling_2014}, travelers often prioritize total travel time and ease of connections, factors that airlines with strategically integrated multimodal networks can leverage.
Overall, the incorporation of airside bus service into airlines' networks reflects a strategic tension between cost efficiency, market coverage, and customer experience. Airlines that successfully navigate these trade-offs ultimately stand to improve their financial sustainability and resilience in the face of regional air service challenges.

\subsection{Passenger Response}
Ultimately, the feasibility of this concept comes down to the passengers' willingness to leverage this service as opposed to waiting for a traditional regional flight or driving directly to their local hub airport. As described at the beginning of this section, this concept has seen real-world applications in North America and Europe. Preliminary reports have shown that in the Northeastern US, the Landline Company has moved a significant number of passengers, demonstrating broad public interest and commercial viability. 

One such example has been on the route established in 2024 between Wilmington Airport (ILG) and Philadelphia International Airport (PHL). Historically, a commercial flight has not existed between these two airports, as they are separated by only 24 miles. Yet, since its establishment, the route has served over 20,000 passengers through December of 2025~\citep{irizarry_shuttle_2026}. 
Inherently, due to their proximity, ILG and PHL reside within overlapping catchment areas. As such, due to PHL's significantly larger connectivity, competitive pricing, and high service frequency, it casts a traffic shadow over the region.
As \citet{fuellhart_airport_2007} discussed, concerning a similar airport pair in the Mid-Atlantic US, this spatial dynamic results in a significant degree of outbound passenger leakage from the regional airport to the hub airport. 
We posit that the observed passenger movements can be attributed to the fact that the small regional airport, ILG, is recapturing a portion of this leakage. This is due to the fact that the advent of this airside-to-airside bus service improves ground access by taking advantage of accessible parking, shorter security lines, and integrated baggage routing at the regional airport~\citep{megan_steckler_new_2023, perez-castells_even_2024}. This reduction in overall travel disutility, in turn, incentivizes passengers to leverage this service rather than flying directly out of PHL~\citep{loo_passengers_2008, ryerson_drive_2018}.

\section{Concluding Remarks} \label{sec:conclusion}
In this work, we presented a mixed integer linear program that aimed to adapt an existing air transportation network such that it may be augmented by airside-to-airside bus operations in lieu of standard regional air service. To validate our model, we conducted a sensitivity analysis to examine the impacts of a variety of variable input parameters and estimated the operational impacts of our proposal through the use of an Agent Based Model Simulation, where significant reductions in average hourly passenger delays were identified in the presence of the augmented network. Lastly, we discussed the current precedents and political hurdles that must be overcome to implement this proposed system.

\subsection{Limitations} \label{ssec:limitations}
Our model and subsequent analysis rely on several parameter assumptions, which inherently limit the broad interpretation of the final results.
In our optimization model, we did not directly consider historical schedules, non-uniform aircraft fleets, and individual airline or policy considerations.
Furthermore, while the current objective function establishes a baseline by minimizing daily operational costs under nominal conditions, it does not presently account for varied operational costs under a diverse set of disruptive scenarios. 
Additionally, as the region under consideration was a reduced portion of the complete air transportation network, it is feasible to assume that the routing of passengers in our optimization result may not reflect the real world, as, in reality, such is impacted by the broader range of flows through the network.
Furthermore, our model simplified demand throughout the fiscal quarter by defining a synthetic representative operational day. 
While this does provide a strong baseline for network planners to utilize, it inherently causes the model to lose fidelity, as this introduces uncertainty with respect to flight schedules and passenger movements that do not occur daily.

With regard to our simulation, the extent to which the excluded flights originating or arriving at the airports under consideration may have ultimately impacted delays beyond those observed in our results is not directly quantified. This is highlighted by the fact that we primarily use true departure time for historical flights to capture potential delays not originating from capacity reductions.
Additionally, although we observed that an overall reduction in average hourly passenger delays occurred through the introduction of bus service in lieu of regional flights, the precise degree and consistency of these improvements may vary throughout the calendar year. In particular, various seasonal variations and demand shifts may increase or reduce the observed result across a broader range of analyzed days.
Lastly, various simplifications were made for efficiency and to serve as a proof of concept for the proposed network. Further investigation is required with regard to the fidelity of our simulation to better understand the complex interactions that airside bus service might produce. 

\subsection{Future Work}
As this research continues, we propose the following suggestions to address the discussed limitations. 
Firstly, the applied estimation process introduces uncertainty not only in the presented results but also in potential downstream applications concerning passenger recovery. 
This matters in particular for routes where small changes in the level of demand can significantly impact their existence. 
Further investigation concerning the reduction of this uncertainty, in conjunction with a quantification of said uncertainty in the final result, is required.
A promising approach to address this is transitioning to a more holistic demand model. This would allow for the characterization of the network while considering consistent demand beyond historical air travel.

Second, while this work incorporated a brief discussion on the notion of passenger willingness to pay for this concept, a direct empirical investigation into these behavioral dynamics was beyond the scope of this study. As this work focuses on a modal transition in the U.S., it is assumed that there is a cultural bias against bus lines due to perceived inconvenience compared to driving or flying. However, travelers often prioritize total travel time and ease of connections, factors that a strategically integrated multimodal network could improve, as discussed in this work. Nonetheless, future behavioral studies and empirical surveys are necessary to quantify how passengers weigh these competing factors and to validate the adoption rates of airside bus services.

Lastly, to enhance the model's practical applicability, future case studies should be expanded to better reflect varied real-world operations and scenarios. This should be incorporated into the sensitivity analysis, which would enable a stronger understanding of the manner in which the proposed network responds to a dynamic environment.
Simultaneously, the optimization model should be expanded to better reflect the air carrier decision-making process and examine how these changes impact not only the observed network topology but also the ultimate resilience of the air transportation network. 
Specifically, future formulations may benefit from the direct coupling of the simulation results to the estimation process. 
By integrating the agent-based simulation and the network optimization, the model could iteratively refine topologies based on disruption-specific operational costs and passenger delays. Furthermore, expanding this framework would provide a natural pathway to relax the initial constraint restricting vehicles to a given O-D pair in the network optimization. This would potentially allow the model to capture the flexible, multi-leg routing realities of an integrated air transportation network.

\clearpage
\appendix
\section{Additional Simulation Visualizations} \label{apdx:addl-viz}

\begin{figure}[pos=htb!]
    \centering
    \includegraphics[width=0.9\linewidth]{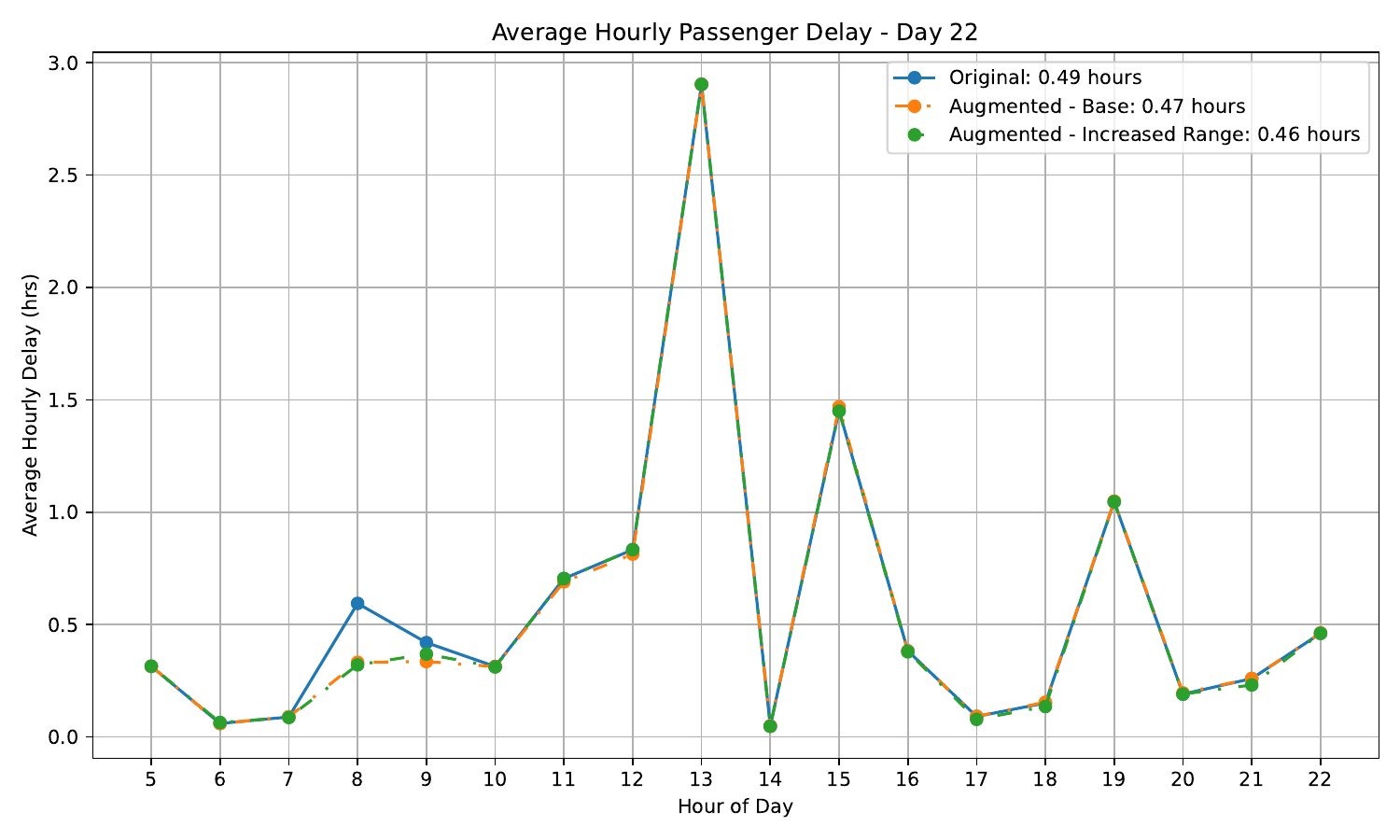}
    \caption{Average hourly passenger delay for both the original and augmented New England network(s) for the nominal case on July 22\textsuperscript{nd}, 2023.}
    \label{fig:sim-results-22}
\end{figure}

\begin{figure}[pos=htb!]
    \centering
    \includegraphics[width=0.9\linewidth]{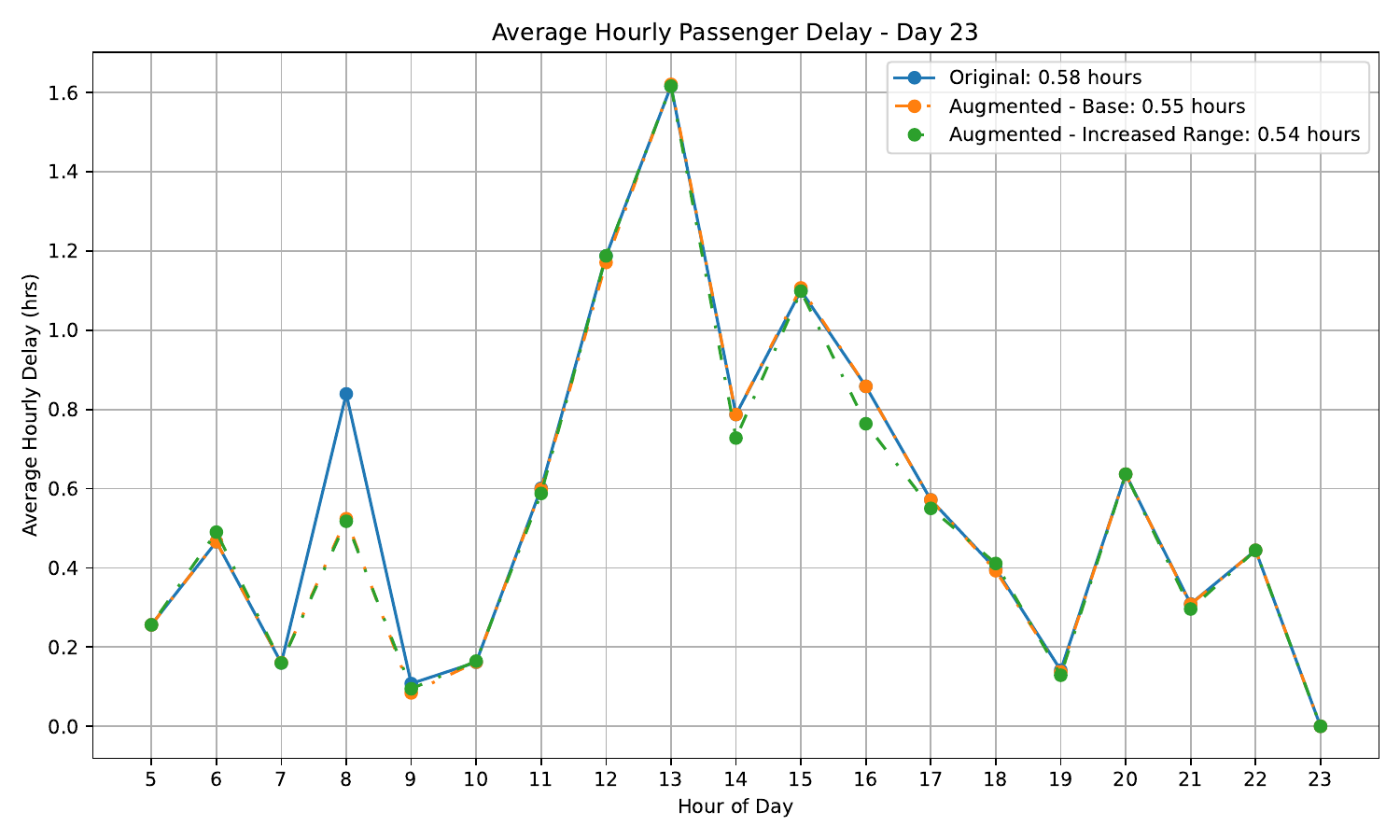}
    \caption{Average hourly passenger delay for both the original and augmented New England network(s) for the nominal case on July 23\textsuperscript{rd}, 2023.}
    \label{fig:sim-results-23}
\end{figure}

\begin{figure}[pos=htb!]
    \centering
    \includegraphics[width=0.9\linewidth]{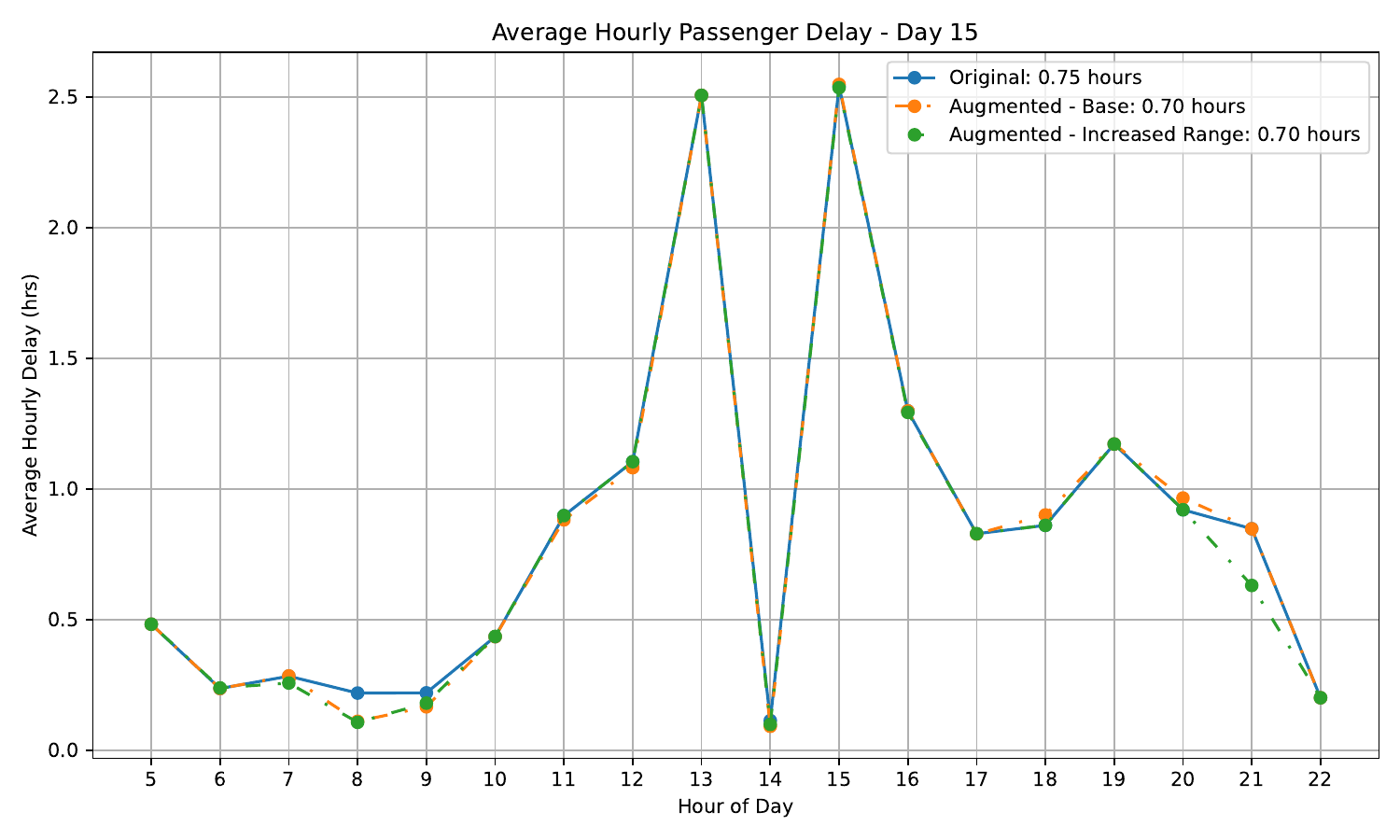}
    \caption{Average hourly passenger delay for both the original and augmented New England network(s) for the nominal case on July 15\textsuperscript{th}, 2023.}
    \label{fig:sim-results-15}
\end{figure}

\begin{figure}[pos=htb!]
    \centering
    \includegraphics[width=0.9\linewidth]{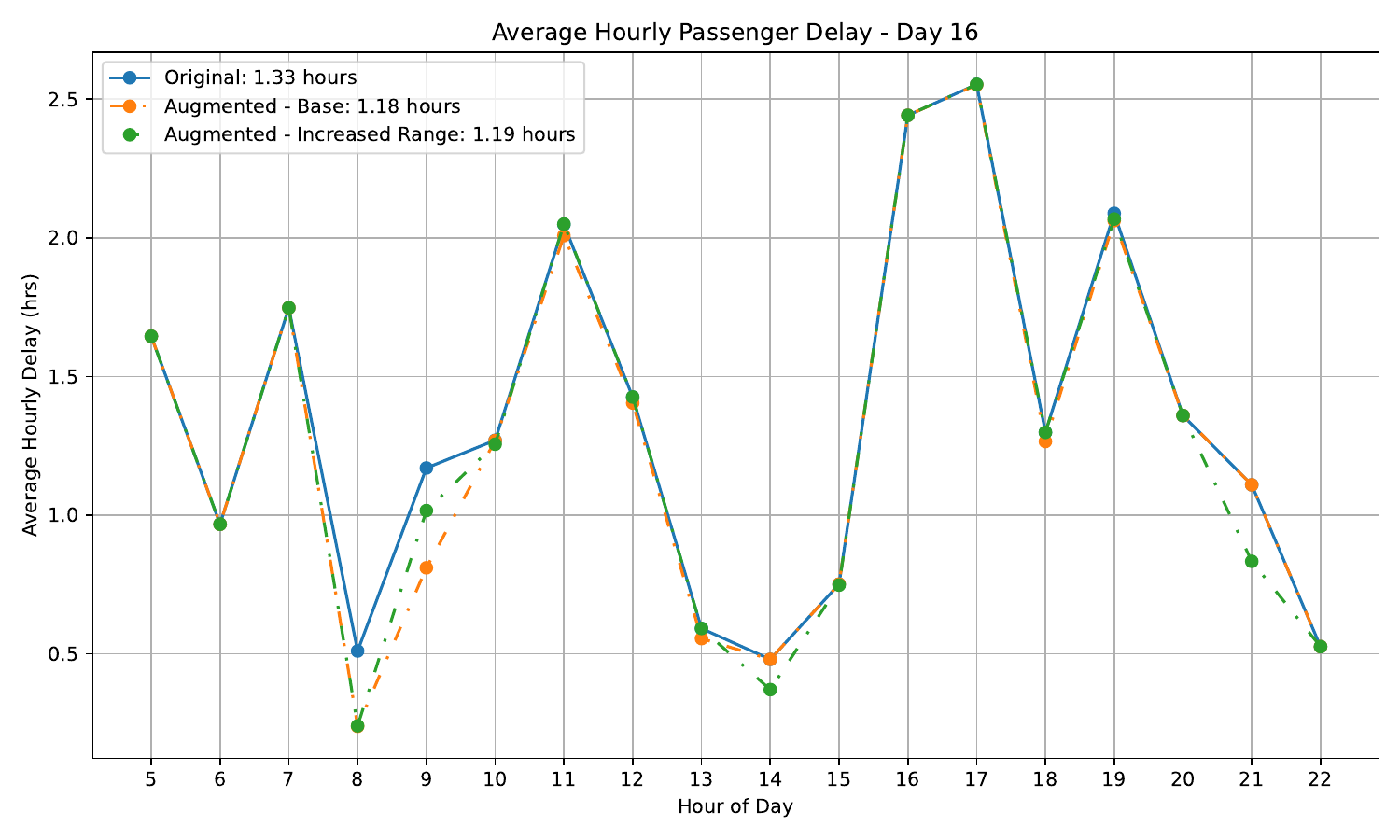}
    \caption{Average hourly passenger delay for both the original and augmented New England network(s) for the nominal case on July 16\textsuperscript{th}, 2023.}
    \label{fig:sim-results-16}
\end{figure}

\clearpage

\printcredits

\section*{Acknowledgments}
The authors thank Jing Xu at the University of California, Berkeley, and the team supporting the NASA University Leadership Initiative (ULI) Center for Air Transportation Resilience (CATRes) in providing relevant data for the Agent Based Model Simulation.

This material is based upon work supported by the National Science Foundation Graduate Research Fellowship Program under Grant No. DGE 2241144 and the National Aeronautics and Space Administration University Leadership initiative under Grant No. 80NSSC24M0068. Any opinions, findings, conclusions or recommendations expressed in this material are those of the author(s) and do not necessarily reflect the views of the National Science Foundation or the National Aeronautics and Space Administration.

\bibliographystyle{cas-model2-names}

\bibliography{references}

\begin{thebibliography}{80}
\expandafter\ifx\csname natexlab\endcsname\relax\def\natexlab#1{#1}\fi
\providecommand{\url}[1]{\texttt{#1}}
\providecommand{\href}[2]{#2}
\providecommand{\path}[1]{#1}
\providecommand{\DOIprefix}{doi:}
\providecommand{\ArXivprefix}{arXiv:}
\providecommand{\URLprefix}{URL: }
\providecommand{\Pubmedprefix}{pmid:}
\providecommand{\doi}[1]{\href{http://dx.doi.org/#1}{\path{#1}}}
\providecommand{\Pubmed}[1]{\href{pmid:#1}{\path{#1}}}
\providecommand{\bibinfo}[2]{#2}
\ifx\xfnm\relax \def\xfnm[#1]{\unskip,\space#1}\fi
\bibitem[{Abdinnour-Helm(1998)}]{abdinnour-helm_hybrid_1998}
\bibinfo{author}{Abdinnour-Helm, S.}, \bibinfo{year}{1998}.
\newblock \bibinfo{title}{A hybrid heuristic for the uncapacitated hub location
  problem}.
\newblock \bibinfo{journal}{European Journal of Operational Research}
  \bibinfo{volume}{106}, \bibinfo{pages}{489--499}.
\newblock \URLprefix
  \url{https://www.sciencedirect.com/science/article/pii/S0377221797002865},
  \DOIprefix\doi{10.1016/S0377-2217(97)00286-5}.
\bibitem[{{Air Canada}(2025)}]{air_canada_air_2025}
\bibinfo{author}{{Air Canada}}, \bibinfo{year}{2025}.
\newblock \bibinfo{title}{Air {Canada} - {Official} {Website}}.
\newblock \URLprefix \url{https://www.aircanada.com/home/ca/en/aco/home.html}.
\bibitem[{{American Airlines}(2025)}]{american_airlines_american_2025}
\bibinfo{author}{{American Airlines}}, \bibinfo{year}{2025}.
\newblock \bibinfo{title}{American {Airlines} - {Airline} tickets and low fares
  at aa.com}.
\newblock \URLprefix \url{https://www.aa.com/homePage.do}.
\bibitem[{{American Bus
  Association}(2025)}]{american_bus_association_motorcoach_2025}
\bibinfo{author}{{American Bus Association}}, \bibinfo{year}{2025}.
\newblock \bibinfo{title}{Motorcoach, {Travel}, \& {Tourism} {Industry}
  {Research}}.
\newblock \URLprefix
  \url{https://www.buses.org/aba-foundation-research-summary/}.
\bibitem[{Archetti et~al.(2022)Archetti, Peirano and
  Speranza}]{archetti_optimization_2022}
\bibinfo{author}{Archetti, C.}, \bibinfo{author}{Peirano, L.},
  \bibinfo{author}{Speranza, M.G.}, \bibinfo{year}{2022}.
\newblock \bibinfo{title}{Optimization in multimodal freight transportation
  problems: {A} {Survey}}.
\newblock \bibinfo{journal}{European Journal of Operational Research}
  \bibinfo{volume}{299}, \bibinfo{pages}{1--20}.
\newblock \URLprefix
  \url{https://www.sciencedirect.com/science/article/pii/S0377221721006263},
  \DOIprefix\doi{10.1016/j.ejor.2021.07.031}.
\bibitem[{Bania et~al.(1998)Bania, Bauer and Zlatoper}]{bania_us_1998}
\bibinfo{author}{Bania, N.}, \bibinfo{author}{Bauer, P.W.},
  \bibinfo{author}{Zlatoper, T.J.}, \bibinfo{year}{1998}.
\newblock \bibinfo{title}{U.{S}. {Air} {Passenger} {Service}: a {Taxonomy} of
  {Route} {Networks}, {Hub} {Locations}, and {Competition}}.
\newblock \bibinfo{journal}{Transportation Research Part E: Logistics and
  Transportation Review} \bibinfo{volume}{34}, \bibinfo{pages}{53--74}.
\newblock \URLprefix
  \url{https://www.sciencedirect.com/science/article/pii/S1366554597000379},
  \DOIprefix\doi{10.1016/S1366-5545(97)00037-9}.
\bibitem[{Barnhart et~al.(2014)Barnhart, Fearing and
  Vaze}]{barnhart_modeling_2014}
\bibinfo{author}{Barnhart, C.}, \bibinfo{author}{Fearing, D.},
  \bibinfo{author}{Vaze, V.}, \bibinfo{year}{2014}.
\newblock \bibinfo{title}{Modeling {Passenger} {Travel} and {Delays} in the
  {National} {Air} {Transportation} {System}}.
\newblock \bibinfo{journal}{Operations Research} \bibinfo{volume}{62},
  \bibinfo{pages}{580--601}.
\newblock \URLprefix
  \url{https://pubsonline.informs.org/doi/10.1287/opre.2014.1268},
  \DOIprefix\doi{10.1287/opre.2014.1268}.
\bibitem[{Bertsimas et~al.(2021)Bertsimas, Ng and
  Yan}]{bertsimas_data-driven_2021}
\bibinfo{author}{Bertsimas, D.}, \bibinfo{author}{Ng, Y.S.},
  \bibinfo{author}{Yan, J.}, \bibinfo{year}{2021}.
\newblock \bibinfo{title}{Data-{Driven} {Transit} {Network} {Design} at
  {Scale}}.
\newblock \bibinfo{journal}{Operations Research} \URLprefix
  \url{https://pubsonline.informs.org/doi/abs/10.1287/opre.2020.2057},
  \DOIprefix\doi{10.1287/opre.2020.2057}.
\bibitem[{{Boeing}(2025)}]{boeing_airplane_2025}
\bibinfo{author}{{Boeing}}, \bibinfo{year}{2025}.
\newblock \bibinfo{title}{Airplane {Characteristics} for {Airport} {Planning}}.
\newblock \URLprefix
  \url{https://www.boeing.com/content/theboeingcompany/us/en/commercial/airports/plan-manuals}.
\bibitem[{Bratu and Barnhart(2005)}]{bratu_analysis_2005}
\bibinfo{author}{Bratu, S.}, \bibinfo{author}{Barnhart, C.},
  \bibinfo{year}{2005}.
\newblock \bibinfo{title}{An analysis of passenger delays using flight
  operations and passenger booking data}.
\newblock \bibinfo{journal}{Air Traffic Control Quarterly}
  \bibinfo{volume}{13}, \bibinfo{pages}{1--27}.
\bibitem[{Brey and Walker(2011)}]{brey_latent_2011}
\bibinfo{author}{Brey, R.}, \bibinfo{author}{Walker, J.L.},
  \bibinfo{year}{2011}.
\newblock \bibinfo{title}{Latent temporal preferences: {An} application to
  airline travel}.
\newblock \bibinfo{journal}{Transportation Research Part A: Policy and
  Practice} \bibinfo{volume}{45}, \bibinfo{pages}{880--895}.
\newblock \URLprefix
  \url{https://www.sciencedirect.com/science/article/pii/S0965856411000681},
  \DOIprefix\doi{10.1016/j.tra.2011.04.010}.
\bibitem[{Bryan and O'Kelly(1999)}]{bryan_hub-and-spoke_1999}
\bibinfo{author}{Bryan, D.L.}, \bibinfo{author}{O'Kelly, M.E.},
  \bibinfo{year}{1999}.
\newblock \bibinfo{title}{Hub-and-{Spoke} {Networks} in {Air} {Transportation}:
  {An} {Analytical} {Review}}.
\newblock \bibinfo{journal}{Journal of Regional Science} \bibinfo{volume}{39},
  \bibinfo{pages}{275--295}.
\newblock \URLprefix
  \url{https://onlinelibrary.wiley.com/doi/abs/10.1111/1467-9787.00134},
  \DOIprefix\doi{10.1111/1467-9787.00134}. \bibinfo{note}{\_eprint:
  https://onlinelibrary.wiley.com/doi/pdf/10.1111/1467-9787.00134}.
\bibitem[{{Bureau of Transportation
  Statistics}(2024)}]{bureau_of_transportation_statistics_number_2024}
\bibinfo{author}{{Bureau of Transportation Statistics}}, \bibinfo{year}{2024}.
\newblock \bibinfo{title}{Number 340 - {Air} {Carrier} {Groups} 2025 {\textbar}
  {Bureau} of {Transportation} {Statistics}}.
\newblock \URLprefix
  \url{https://www.bts.gov/explore-topics-and-geography/modes/aviation/number-340-air-carrier-groups-2025}.
\bibitem[{{Bureau of Transportation
  Statistics}(2025)}]{bureau_of_transportation_statistics_transstats_2025}
\bibinfo{author}{{Bureau of Transportation Statistics}}, \bibinfo{year}{2025}.
\newblock \bibinfo{title}{{TransStats}}.
\newblock \URLprefix \url{https://www.transtats.bts.gov/}.
\bibitem[{{Bureau of Transportation
  Statistics}(2026)}]{bureau_of_transportation_statistics_origin-destination_2026}
\bibinfo{author}{{Bureau of Transportation Statistics}}, \bibinfo{year}{2026}.
\newblock \bibinfo{title}{Origin-{Destination} 40 {Percent} {Survey} of
  {Airline} {Passenger} {Traffic} {\textbar} {Bureau} of {Transportation}
  {Statistics}}.
\newblock \URLprefix \url{https://www.bts.gov/OD-40}.
\bibitem[{Campbell(1994)}]{campbell_integer_1994}
\bibinfo{author}{Campbell, J.F.}, \bibinfo{year}{1994}.
\newblock \bibinfo{title}{Integer programming formulations of discrete hub
  location problems}.
\newblock \bibinfo{journal}{European Journal of Operational Research}
  \bibinfo{volume}{72}, \bibinfo{pages}{387--405}.
\newblock \URLprefix
  \url{https://www.sciencedirect.com/science/article/pii/0377221794903182},
  \DOIprefix\doi{10.1016/0377-2217(94)90318-2}.
\bibitem[{Campbell and O'Kelly(2012)}]{campbell_twenty-five_2012}
\bibinfo{author}{Campbell, J.F.}, \bibinfo{author}{O'Kelly, M.E.},
  \bibinfo{year}{2012}.
\newblock \bibinfo{title}{Twenty-{Five} {Years} of {Hub} {Location}
  {Research}}.
\newblock \bibinfo{journal}{Transportation Science} \bibinfo{volume}{46},
  \bibinfo{pages}{153--169}.
\newblock \URLprefix
  \url{https://pubsonline.informs.org/doi/10.1287/trsc.1120.0410},
  \DOIprefix\doi{10.1287/trsc.1120.0410}.
\bibitem[{Chau and Gkiotsalitis(2025)}]{chau_systematic_2025}
\bibinfo{author}{Chau, M.L.Y.}, \bibinfo{author}{Gkiotsalitis, K.},
  \bibinfo{year}{2025}.
\newblock \bibinfo{title}{A systematic literature review on the use of
  metaheuristics for the optimisation of multimodal transportation}.
\newblock \bibinfo{journal}{Evolutionary Intelligence} \bibinfo{volume}{18},
  \bibinfo{pages}{36}.
\newblock \URLprefix \url{https://doi.org/10.1007/s12065-025-01020-2},
  \DOIprefix\doi{10.1007/s12065-025-01020-2}.
\bibitem[{Clausen et~al.(2010)Clausen, Larsen, Larsen and
  Rezanova}]{clausen_disruption_2010}
\bibinfo{author}{Clausen, J.}, \bibinfo{author}{Larsen, A.},
  \bibinfo{author}{Larsen, J.}, \bibinfo{author}{Rezanova, N.J.},
  \bibinfo{year}{2010}.
\newblock \bibinfo{title}{Disruption management in the airline
  industry—{Concepts}, models and methods}.
\newblock \bibinfo{journal}{Computers \& Operations Research}
  \bibinfo{volume}{37}, \bibinfo{pages}{809--821}.
\newblock \URLprefix
  \url{https://www.sciencedirect.com/science/article/pii/S0305054809000914},
  \DOIprefix\doi{10.1016/j.cor.2009.03.027}.
\bibitem[{Contreras et~al.(2011a)Contreras, Cordeau and
  Laporte}]{contreras_stochastic_2011}
\bibinfo{author}{Contreras, I.}, \bibinfo{author}{Cordeau, J.F.},
  \bibinfo{author}{Laporte, G.}, \bibinfo{year}{2011}a.
\newblock \bibinfo{title}{Stochastic uncapacitated hub location}.
\newblock \bibinfo{journal}{European Journal of Operational Research}
  \bibinfo{volume}{212}, \bibinfo{pages}{518--528}.
\newblock \URLprefix
  \url{https://www.sciencedirect.com/science/article/pii/S0377221711001494},
  \DOIprefix\doi{10.1016/j.ejor.2011.02.018}.
\bibitem[{Contreras et~al.(2011b)Contreras, Díaz and
  Fernández}]{contreras_branch_2011}
\bibinfo{author}{Contreras, I.}, \bibinfo{author}{Díaz, J.A.},
  \bibinfo{author}{Fernández, E.}, \bibinfo{year}{2011}b.
\newblock \bibinfo{title}{Branch and {Price} for {Large}-{Scale} {Capacitated}
  {Hub} {Location} {Problems} with {Single} {Assignment}}.
\newblock \bibinfo{journal}{INFORMS Journal on Computing} \bibinfo{volume}{23},
  \bibinfo{pages}{41--55}.
\newblock \URLprefix
  \url{https://pubsonline.informs.org/doi/abs/10.1287/ijoc.1100.0391},
  \DOIprefix\doi{10.1287/ijoc.1100.0391}.
\bibitem[{Delgado et~al.(2026)Delgado, Gurtner, Weiszer, Bolić and
  Cook}]{delgado_mercury_2026}
\bibinfo{author}{Delgado, L.}, \bibinfo{author}{Gurtner, G.},
  \bibinfo{author}{Weiszer, M.}, \bibinfo{author}{Bolić, T.},
  \bibinfo{author}{Cook, A.}, \bibinfo{year}{2026}.
\newblock \bibinfo{title}{Mercury: {An} open-source simulator for the
  evaluation of air transport mobility}.
\newblock \bibinfo{journal}{European Journal of Transport and Infrastructure
  Research} \bibinfo{volume}{26}.
\newblock \URLprefix
  \url{https://journals.open.tudelft.nl/ejtir/article/view/7529},
  \DOIprefix\doi{10.59490/ejtir.2026.26.1.7529}.
\bibitem[{{Delta Air Lines}(2024)}]{delta_air_lines_sec_2024}
\bibinfo{author}{{Delta Air Lines}}, \bibinfo{year}{2024}.
\newblock \bibinfo{title}{{SEC} {Form} 8-{K}}.
\newblock \bibinfo{type}{Technical Report}.
\newblock \URLprefix
  \url{https://www.sec.gov/Archives/edgar/data/27904/000168316824005369/delta_8k.htm}.
\bibitem[{{Department of Transportation, Office of the
  Secretary}(2023)}]{department_of_transportation_office_of_the_secretary_updates_2023}
\bibinfo{author}{{Department of Transportation, Office of the Secretary}},
  \bibinfo{year}{2023}.
\newblock \bibinfo{title}{Updates to the {Origin}-{Destination} {Survey} of
  {Airline} {Passengers}}.
\newblock \bibinfo{type}{Technical Report} \bibinfo{number}{FR Doc.
  2022-28535}.
\newblock \URLprefix
  \url{https://www.federalregister.gov/documents/2023/01/31/2022-28535/updates-to-the-origin-destination-survey-of-airline-passengers}.
  \bibinfo{note}{final Rule; 88 Fed. Reg. 6145; to be codified at 14 C.F.R.
  pts. 241 \& 298}.
\bibitem[{Durán-Micco and Vansteenwegen(2022)}]{duran-micco_survey_2022}
\bibinfo{author}{Durán-Micco, J.}, \bibinfo{author}{Vansteenwegen, P.},
  \bibinfo{year}{2022}.
\newblock \bibinfo{title}{A survey on the transit network design and frequency
  setting problem}.
\newblock \bibinfo{journal}{Public Transport} \bibinfo{volume}{14},
  \bibinfo{pages}{155--190}.
\newblock \URLprefix \url{https://doi.org/10.1007/s12469-021-00284-y},
  \DOIprefix\doi{10.1007/s12469-021-00284-y}.
\bibitem[{Ennen et~al.(2019)Ennen, Allroggen and Malina}]{ennen_non-stop_2019}
\bibinfo{author}{Ennen, D.}, \bibinfo{author}{Allroggen, F.},
  \bibinfo{author}{Malina, R.}, \bibinfo{year}{2019}.
\newblock \bibinfo{title}{Non-stop versus connecting air services: {Airfares},
  costs, and consumers’ willingness to pay} \URLprefix
  \url{https://dspace.mit.edu/handle/1721.1/121459}. \bibinfo{note}{accepted:
  2019-06-28T19:49:56Z}.
\bibitem[{{FAA}(2022)}]{faa_airport_2022}
\bibinfo{author}{{FAA}}, \bibinfo{year}{2022}.
\newblock \bibinfo{title}{Airport {Categories}}.
\newblock \URLprefix
  \url{https://www.faa.gov/airports/planning_capacity/categories}.
\bibitem[{FAA(2024)}]{faa_faa_2024}
\bibinfo{author}{FAA}, \bibinfo{year}{2024}.
\newblock \bibinfo{title}{{FAA} {Order} {JO} 7210.{3DD} - {Facility}
  {Operation} and {Administration}}.
\newblock \URLprefix
  \url{https://www.faa.gov/air_traffic/publications/atpubs/foa_html/chap18_section_7.html}.
\bibitem[{{FAA}(2025a)}]{faa_aircraft_2025}
\bibinfo{author}{{FAA}}, \bibinfo{year}{2025}a.
\newblock \bibinfo{title}{Aircraft {Registration} {Database}}.
\newblock \URLprefix \url{https://registry.faa.gov/aircraftinquiry}.
\bibitem[{{FAA}(2025b)}]{faa_aspm_2025}
\bibinfo{author}{{FAA}}, \bibinfo{year}{2025}b.
\newblock \bibinfo{title}{{ASPM} {Taxi} {Times}: {Standard} {Report} -
  {ASPMHelp}}.
\newblock \URLprefix
  \url{https://www.aspm.faa.gov/aspmhelp/index/ASPM_Taxi_Times__Standard_Report.html}.
\bibitem[{Farahani et~al.(2013a)Farahani, Hekmatfar, Arabani and
  Nikbakhsh}]{farahani_hub_2013}
\bibinfo{author}{Farahani, R.Z.}, \bibinfo{author}{Hekmatfar, M.},
  \bibinfo{author}{Arabani, A.B.}, \bibinfo{author}{Nikbakhsh, E.},
  \bibinfo{year}{2013}a.
\newblock \bibinfo{title}{Hub location problems: {A} review of models,
  classification, solution techniques, and applications}.
\newblock \bibinfo{journal}{Computers \& Industrial Engineering}
  \bibinfo{volume}{64}, \bibinfo{pages}{1096--1109}.
\newblock \URLprefix
  \url{https://www.sciencedirect.com/science/article/pii/S0360835213000326},
  \DOIprefix\doi{10.1016/j.cie.2013.01.012}.
\bibitem[{Farahani et~al.(2013b)Farahani, Miandoabchi, Szeto and
  Rashidi}]{farahani_review_2013}
\bibinfo{author}{Farahani, R.Z.}, \bibinfo{author}{Miandoabchi, E.},
  \bibinfo{author}{Szeto, W.Y.}, \bibinfo{author}{Rashidi, H.},
  \bibinfo{year}{2013}b.
\newblock \bibinfo{title}{A review of urban transportation network design
  problems}.
\newblock \bibinfo{journal}{European Journal of Operational Research}
  \bibinfo{volume}{229}, \bibinfo{pages}{281--302}.
\newblock \URLprefix
  \url{https://www.sciencedirect.com/science/article/pii/S0377221713000106},
  \DOIprefix\doi{10.1016/j.ejor.2013.01.001}.
\bibitem[{Fuellhart(2007)}]{fuellhart_airport_2007}
\bibinfo{author}{Fuellhart, K.}, \bibinfo{year}{2007}.
\newblock \bibinfo{title}{Airport catchment and leakage in a multi-airport
  region: {The} case of {Harrisburg} {International}}.
\newblock \bibinfo{journal}{Journal of Transport Geography}
  \bibinfo{volume}{15}, \bibinfo{pages}{231--244}.
\newblock \URLprefix
  \url{https://www.sciencedirect.com/science/article/pii/S096669230600072X},
  \DOIprefix\doi{10.1016/j.jtrangeo.2006.08.001}.
\bibitem[{Gao(2021)}]{gao_what_2021}
\bibinfo{author}{Gao, Y.}, \bibinfo{year}{2021}.
\newblock \bibinfo{title}{What is the busiest time at an airport? {Clustering}
  {U}.{S}. hub airports based on passenger movements}.
\newblock \bibinfo{journal}{Journal of Transport Geography}
  \bibinfo{volume}{90}, \bibinfo{pages}{102931}.
\newblock \URLprefix
  \url{https://www.sciencedirect.com/science/article/pii/S0966692320310085},
  \DOIprefix\doi{10.1016/j.jtrangeo.2020.102931}.
\bibitem[{Gertsbakh and Shpungin(2011)}]{gertsbakh_network_2011}
\bibinfo{author}{Gertsbakh, I.}, \bibinfo{author}{Shpungin, Y.},
  \bibinfo{year}{2011}.
\newblock \bibinfo{title}{Network {Reliability} and {Resilience}}.
\newblock \bibinfo{publisher}{Springer}, \bibinfo{address}{Berlin, Heidelberg}.
\newblock \URLprefix \url{https://link.springer.com/10.1007/978-3-642-22374-7},
  \DOIprefix\doi{10.1007/978-3-642-22374-7}.
\bibitem[{Gurtner et~al.(2021)Gurtner, Delgado and
  Valput}]{gurtner_agent-based_2021}
\bibinfo{author}{Gurtner, G.}, \bibinfo{author}{Delgado, L.},
  \bibinfo{author}{Valput, D.}, \bibinfo{year}{2021}.
\newblock \bibinfo{title}{An agent-based model for air transportation to
  capture network effects in assessing delay management mechanisms}.
\newblock \bibinfo{journal}{Transportation Research Part C: Emerging
  Technologies} \bibinfo{volume}{133}, \bibinfo{pages}{103358}.
\newblock \URLprefix
  \url{https://www.sciencedirect.com/science/article/pii/S0968090X21003600},
  \DOIprefix\doi{10.1016/j.trc.2021.103358}.
\bibitem[{Hartley(2025)}]{hartley_american_2025}
\bibinfo{author}{Hartley, P.}, \bibinfo{year}{2025}.
\newblock \bibinfo{title}{American {Airlines} '{Landline}' {Service} {Gets}
  {Passengers} {To} \& {From} {Philadelphia} {Via} {Bus}}.
\newblock \URLprefix
  \url{https://simpleflying.com/american-airlines-landline-philadelphia-bus-service/}.
  \bibinfo{note}{section: Airlines}.
\bibitem[{Hosseini et~al.(2016)Hosseini, Barker and
  Ramirez-Marquez}]{hosseini_review_2016}
\bibinfo{author}{Hosseini, S.}, \bibinfo{author}{Barker, K.},
  \bibinfo{author}{Ramirez-Marquez, J.E.}, \bibinfo{year}{2016}.
\newblock \bibinfo{title}{A review of definitions and measures of system
  resilience}.
\newblock \bibinfo{journal}{Reliability Engineering \& System Safety}
  \bibinfo{volume}{145}, \bibinfo{pages}{47--61}.
\newblock \URLprefix
  \url{https://www.sciencedirect.com/science/article/pii/S0951832015002483},
  \DOIprefix\doi{10.1016/j.ress.2015.08.006}.
\bibitem[{Huang et~al.(2022)Huang, Cui, Zhang, Tong, Shi and
  Liu}]{huang_overview_2022}
\bibinfo{author}{Huang, J.}, \bibinfo{author}{Cui, Y.}, \bibinfo{author}{Zhang,
  L.}, \bibinfo{author}{Tong, W.}, \bibinfo{author}{Shi, Y.},
  \bibinfo{author}{Liu, Z.}, \bibinfo{year}{2022}.
\newblock \bibinfo{title}{An {Overview} of {Agent}-{Based} {Models} for
  {Transport} {Simulation} and {Analysis}}.
\newblock \bibinfo{journal}{Journal of Advanced Transportation}
  \bibinfo{volume}{2022}, \bibinfo{pages}{1252534}.
\newblock \URLprefix
  \url{https://onlinelibrary.wiley.com/doi/abs/10.1155/2022/1252534},
  \DOIprefix\doi{10.1155/2022/1252534}. \bibinfo{note}{\_eprint:
  https://onlinelibrary.wiley.com/doi/pdf/10.1155/2022/1252534}.
\bibitem[{{International Air Transport
  Association}(2025)}]{international_air_transport_association_economics_2025}
\bibinfo{author}{{International Air Transport Association}},
  \bibinfo{year}{2025}.
\newblock \bibinfo{title}{Economics {Report} {Library}}.
\newblock \URLprefix
  \url{https://www.iata.org/en/publications/economics/economics-library/}.
\bibitem[{Irizarry(2026)}]{irizarry_shuttle_2026}
\bibinfo{author}{Irizarry, J.}, \bibinfo{year}{2026}.
\newblock \bibinfo{title}{The shuttle service taking passengers to
  {Philadelphia} {International} {Airport} hits 20,000 customers}.
\newblock \URLprefix
  \url{https://www.delawarepublic.org/business/2026-01-15/the-shuttle-service-taking-passengers-to-philadelphia-international-airport-hits-20-000-customers}.
  \bibinfo{note}{section: Business}.
\bibitem[{Jacquillat(2022)}]{jacquillat_predictive_2022}
\bibinfo{author}{Jacquillat, A.}, \bibinfo{year}{2022}.
\newblock \bibinfo{title}{Predictive and {Prescriptive} {Analytics} {Toward}
  {Passenger}-{Centric} {Ground} {Delay} {Programs}}.
\newblock \bibinfo{journal}{Transportation Science} \bibinfo{volume}{56},
  \bibinfo{pages}{265--298}.
\newblock \URLprefix
  \url{https://pubsonline.informs.org/doi/abs/10.1287/trsc.2021.1081},
  \DOIprefix\doi{10.1287/trsc.2021.1081}.
\bibitem[{Jaillet et~al.(1996)Jaillet, Song and Yu}]{jaillet_airline_1996}
\bibinfo{author}{Jaillet, P.}, \bibinfo{author}{Song, G.}, \bibinfo{author}{Yu,
  G.}, \bibinfo{year}{1996}.
\newblock \bibinfo{title}{Airline network design and hub location problems}.
\newblock \bibinfo{journal}{Location Science} \bibinfo{volume}{4},
  \bibinfo{pages}{195--212}.
\newblock \URLprefix
  \url{https://www.sciencedirect.com/science/article/pii/S0966834996000162},
  \DOIprefix\doi{10.1016/S0966-8349(96)00016-2}.
\bibitem[{Janić(2005)}]{janic_modeling_2005}
\bibinfo{author}{Janić, M.}, \bibinfo{year}{2005}.
\newblock \bibinfo{title}{Modeling the {Large} {Scale} {Disruptions} of an
  {Airline} {Network}}.
\newblock \bibinfo{journal}{Journal of Transportation Engineering}
  \bibinfo{volume}{131}, \bibinfo{pages}{249--260}.
\newblock \URLprefix
  \url{https://ascelibrary.org/doi/10.1061/%28ASCE%290733-947X%282005%29131%3A4%28249%29},
  \DOIprefix\doi{10.1061/(ASCE)0733-947X(2005)131:4(249)}.
\bibitem[{Janić(2015)}]{janic_reprint_2015}
\bibinfo{author}{Janić, M.}, \bibinfo{year}{2015}.
\newblock \bibinfo{title}{Reprint of “{Modelling} the resilience, friability
  and costs of an air transport network affected by a large-scale disruptive
  event”}.
\newblock \bibinfo{journal}{Transportation Research Part A: Policy and
  Practice} \bibinfo{volume}{81}, \bibinfo{pages}{77--92}.
\newblock \URLprefix
  \url{https://www.sciencedirect.com/science/article/pii/S0965856415002049},
  \DOIprefix\doi{10.1016/j.tra.2015.07.012}.
\bibitem[{Jenelius(2009)}]{jenelius_network_2009}
\bibinfo{author}{Jenelius, E.}, \bibinfo{year}{2009}.
\newblock \bibinfo{title}{Network structure and travel patterns: explaining the
  geographical disparities of road network vulnerability}.
\newblock \bibinfo{journal}{Journal of Transport Geography}
  \bibinfo{volume}{17}, \bibinfo{pages}{234--244}.
\newblock \URLprefix
  \url{https://www.sciencedirect.com/science/article/pii/S0966692308000550},
  \DOIprefix\doi{10.1016/j.jtrangeo.2008.06.002}.
\bibitem[{Jenelius(2010)}]{jenelius_user_2010}
\bibinfo{author}{Jenelius, E.}, \bibinfo{year}{2010}.
\newblock \bibinfo{title}{User inequity implications of road network
  vulnerability}.
\newblock \bibinfo{journal}{Journal of Transport and Land Use}
  \bibinfo{volume}{2}, \bibinfo{pages}{57--73}.
\bibitem[{Jenelius et~al.(2006)Jenelius, Petersen and
  Mattsson}]{jenelius_importance_2006}
\bibinfo{author}{Jenelius, E.}, \bibinfo{author}{Petersen, T.},
  \bibinfo{author}{Mattsson, L.G.}, \bibinfo{year}{2006}.
\newblock \bibinfo{title}{Importance and exposure in road network vulnerability
  analysis}.
\newblock \bibinfo{journal}{Transportation Research Part A: Policy and
  Practice} \bibinfo{volume}{40}, \bibinfo{pages}{537--560}.
\newblock \URLprefix
  \url{https://www.sciencedirect.com/science/article/pii/S096585640500162X},
  \DOIprefix\doi{10.1016/j.tra.2005.11.003}.
\bibitem[{{KLM Royal Dutch
  Airlines}(2025)}]{klm_royal_dutch_airlines_travel_2025}
\bibinfo{author}{{KLM Royal Dutch Airlines}}, \bibinfo{year}{2025}.
\newblock \bibinfo{title}{Travel free to or from {Schiphol} with the {KLM}
  {Bus}}.
\newblock \URLprefix \url{https://bus.klm.nl/en}.
\bibitem[{{Landline Company}(2025)}]{landline_company_landline_2025}
\bibinfo{author}{{Landline Company}}, \bibinfo{year}{2025}.
\newblock \bibinfo{title}{Landline {Partners}}.
\newblock \URLprefix \url{https://landlineco.com/partners/}.
\bibitem[{Li et~al.(2022)Li, Gopalakrishnan, Balakrishnan, Shin, Jalan, Nandi
  and Marla}]{li_dynamics_2022}
\bibinfo{author}{Li, M.Z.}, \bibinfo{author}{Gopalakrishnan, K.},
  \bibinfo{author}{Balakrishnan, H.}, \bibinfo{author}{Shin, S.H.},
  \bibinfo{author}{Jalan, D.}, \bibinfo{author}{Nandi, A.},
  \bibinfo{author}{Marla, L.}, \bibinfo{year}{2022}.
\newblock \bibinfo{title}{Dynamics of disruption and recovery in air
  transportation networks}.
\newblock \bibinfo{journal}{CEAS Aeronautical Journal} \bibinfo{volume}{13},
  \bibinfo{pages}{347--357}.
\newblock \URLprefix \url{https://doi.org/10.1007/s13272-021-00521-x},
  \DOIprefix\doi{10.1007/s13272-021-00521-x}.
\bibitem[{Loo(2008)}]{loo_passengers_2008}
\bibinfo{author}{Loo, B.P.Y.}, \bibinfo{year}{2008}.
\newblock \bibinfo{title}{Passengers’ airport choice within multi-airport
  regions ({MARs}): some insights from a stated preference survey at {Hong}
  {Kong} {International} {Airport}}.
\newblock \bibinfo{journal}{Journal of Transport Geography}
  \bibinfo{volume}{16}, \bibinfo{pages}{117--125}.
\newblock \URLprefix
  \url{https://www.sciencedirect.com/science/article/pii/S0966692307000658},
  \DOIprefix\doi{10.1016/j.jtrangeo.2007.05.003}.
\bibitem[{{Lufthansa}(2025)}]{lufthansa_lufthansa_2025}
\bibinfo{author}{{Lufthansa}}, \bibinfo{year}{2025}.
\newblock \bibinfo{title}{Lufthansa {Express} {Bus}}.
\newblock \URLprefix
  \url{https://www.lufthansa.com/us/en/lufthansa-express-bus}.
\bibitem[{Marla et~al.(2017)Marla, Vaaben and Barnhart}]{marla_integrated_2017}
\bibinfo{author}{Marla, L.}, \bibinfo{author}{Vaaben, B.},
  \bibinfo{author}{Barnhart, C.}, \bibinfo{year}{2017}.
\newblock \bibinfo{title}{Integrated {Disruption} {Management} and {Flight}
  {Planning} to {Trade} {Off} {Delays} and {Fuel} {Burn}}.
\newblock \bibinfo{journal}{Transportation Science} \bibinfo{volume}{51},
  \bibinfo{pages}{88--111}.
\newblock \URLprefix
  \url{https://pubsonline-informs-org.proxy.lib.umich.edu/doi/abs/10.1287/trsc.2015.0609},
  \DOIprefix\doi{10.1287/trsc.2015.0609}.
\bibitem[{Marzuoli et~al.(2016)Marzuoli, Boidot, Colomar, Guerpillon, Feron,
  Bayen and Hansen}]{marzuoli_improving_2016}
\bibinfo{author}{Marzuoli, A.}, \bibinfo{author}{Boidot, E.},
  \bibinfo{author}{Colomar, P.}, \bibinfo{author}{Guerpillon, M.},
  \bibinfo{author}{Feron, E.}, \bibinfo{author}{Bayen, A.},
  \bibinfo{author}{Hansen, M.}, \bibinfo{year}{2016}.
\newblock \bibinfo{title}{Improving {Disruption} {Management} {With}
  {Multimodal} {Collaborative} {Decision}-{Making}: {A} {Case} {Study} of the
  {Asiana} {Crash} and {Lessons} {Learned}}.
\newblock \bibinfo{journal}{IEEE Transactions on Intelligent Transportation
  Systems} \bibinfo{volume}{17}, \bibinfo{pages}{2699--2717}.
\newblock \URLprefix
  \url{https://ieeexplore.ieee.org/abstract/document/7457252},
  \DOIprefix\doi{10.1109/TITS.2016.2536733}.
\bibitem[{Mattsson and Jenelius(2015)}]{mattsson_vulnerability_2015}
\bibinfo{author}{Mattsson, L.G.}, \bibinfo{author}{Jenelius, E.},
  \bibinfo{year}{2015}.
\newblock \bibinfo{title}{Vulnerability and resilience of transport systems –
  {A} discussion of recent research}.
\newblock \bibinfo{journal}{Transportation Research Part A: Policy and
  Practice} \bibinfo{volume}{81}, \bibinfo{pages}{16--34}.
\newblock \URLprefix
  \url{https://www.sciencedirect.com/science/article/pii/S0965856415001603},
  \DOIprefix\doi{10.1016/j.tra.2015.06.002}.
\bibitem[{{Megan Steckler}(2023)}]{megan_steckler_new_2023}
\bibinfo{author}{{Megan Steckler}}, \bibinfo{year}{2023}.
\newblock \bibinfo{title}{New {Security} {Process} {Makes} {American}
  {Airlines}-{Landline} {Travel} {Seamless} {\textbar} {PHL}.org}.
\newblock \URLprefix \url{https://www.phl.org/newsroom/Landline-update}.
\bibitem[{Miller-Hooks et~al.(2012)Miller-Hooks, Zhang and
  Faturechi}]{miller-hooks_measuring_2012}
\bibinfo{author}{Miller-Hooks, E.}, \bibinfo{author}{Zhang, X.},
  \bibinfo{author}{Faturechi, R.}, \bibinfo{year}{2012}.
\newblock \bibinfo{title}{Measuring and maximizing resilience of freight
  transportation networks}.
\newblock \bibinfo{journal}{Computers \& Operations Research}
  \bibinfo{volume}{39}, \bibinfo{pages}{1633--1643}.
\newblock \URLprefix
  \url{https://www.sciencedirect.com/science/article/pii/S0305054811002784},
  \DOIprefix\doi{10.1016/j.cor.2011.09.017}.
\bibitem[{O'Kelly and Miller(1994)}]{okelly_hub_1994}
\bibinfo{author}{O'Kelly, M.E.}, \bibinfo{author}{Miller, H.J.},
  \bibinfo{year}{1994}.
\newblock \bibinfo{title}{The hub network design problem: {A} review and
  synthesis}.
\newblock \bibinfo{journal}{Journal of Transport Geography}
  \bibinfo{volume}{2}, \bibinfo{pages}{31--40}.
\newblock \URLprefix
  \url{https://www.sciencedirect.com/science/article/pii/0966692394900329},
  \DOIprefix\doi{10.1016/0966-6923(94)90032-9}.
\bibitem[{Pauwels et~al.(2024)Pauwels, Buyle and
  Dewulf}]{pauwels_regional_2024}
\bibinfo{author}{Pauwels, J.}, \bibinfo{author}{Buyle, S.},
  \bibinfo{author}{Dewulf, W.}, \bibinfo{year}{2024}.
\newblock \bibinfo{title}{Regional airports revisited: {Unveiling} pressing
  research gaps and proposing a uniform definition}.
\newblock \bibinfo{journal}{Journal of the Air Transport Research Society}
  \bibinfo{volume}{2}, \bibinfo{pages}{100008}.
\newblock \URLprefix
  \url{https://www.sciencedirect.com/science/article/pii/S2941198X24000034},
  \DOIprefix\doi{10.1016/j.jatrs.2024.100008}.
\bibitem[{Perez-Castells(2024)}]{perez-castells_even_2024}
\bibinfo{author}{Perez-Castells, A.}, \bibinfo{year}{2024}.
\newblock \bibinfo{title}{Even more travelers flying {American} {Airlines} will
  be able to skip security lines at {PHL} soon}.
\newblock \URLprefix
  \url{https://www.inquirer.com/business/american-airlines-landline-security-bus-airport-20240518.html}.
  \bibinfo{note}{section: Business, business, business}.
\bibitem[{Rodríguez-Déniz et~al.(2013)Rodríguez-Déniz, Suau-Sanchez and
  Voltes-Dorta}]{rodriguez-deniz_classifying_2013}
\bibinfo{author}{Rodríguez-Déniz, H.}, \bibinfo{author}{Suau-Sanchez, P.},
  \bibinfo{author}{Voltes-Dorta, A.}, \bibinfo{year}{2013}.
\newblock \bibinfo{title}{Classifying airports according to their hub
  dimensions: an application to the {US} domestic network}.
\newblock \bibinfo{journal}{Journal of Transport Geography}
  \bibinfo{volume}{33}, \bibinfo{pages}{188--195}.
\newblock \URLprefix
  \url{https://www.sciencedirect.com/science/article/pii/S0966692313002093},
  \DOIprefix\doi{10.1016/j.jtrangeo.2013.10.011}.
\bibitem[{Rodríguez-Núñez and
  García-Palomares(2014)}]{rodriguez-nunez_measuring_2014}
\bibinfo{author}{Rodríguez-Núñez, E.}, \bibinfo{author}{García-Palomares,
  J.C.}, \bibinfo{year}{2014}.
\newblock \bibinfo{title}{Measuring the vulnerability of public transport
  networks}.
\newblock \bibinfo{journal}{Journal of Transport Geography}
  \bibinfo{volume}{35}, \bibinfo{pages}{50--63}.
\newblock \URLprefix
  \url{https://www.sciencedirect.com/science/article/pii/S0966692314000180},
  \DOIprefix\doi{10.1016/j.jtrangeo.2014.01.008}.
\bibitem[{Ryerson and Kim(2018)}]{ryerson_drive_2018}
\bibinfo{author}{Ryerson, M.S.}, \bibinfo{author}{Kim, A.M.},
  \bibinfo{year}{2018}.
\newblock \bibinfo{title}{A drive for better air service: {How} air service
  imbalances across neighboring regions integrate air and highway demands}.
\newblock \bibinfo{journal}{Transportation Research Part A: Policy and
  Practice} \bibinfo{volume}{114}, \bibinfo{pages}{237--255}.
\newblock \URLprefix
  \url{https://www.sciencedirect.com/science/article/pii/S0965856417304573},
  \DOIprefix\doi{10.1016/j.tra.2017.10.005}.
\bibitem[{Ryerson and Kim(2013)}]{ryerson_integrating_2013}
\bibinfo{author}{Ryerson, M.S.}, \bibinfo{author}{Kim, H.},
  \bibinfo{year}{2013}.
\newblock \bibinfo{title}{Integrating airline operational practices into
  passenger airline hub definition}.
\newblock \bibinfo{journal}{Journal of Transport Geography}
  \bibinfo{volume}{31}, \bibinfo{pages}{84--93}.
\newblock \URLprefix
  \url{https://www.sciencedirect.com/science/article/pii/S0966692313001099},
  \DOIprefix\doi{10.1016/j.jtrangeo.2013.05.013}.
\bibitem[{Schäfer et~al.(2014)Schäfer, Strohmeier, Lenders, Martinovic and
  Wilhelm}]{schafer_bringing_2014}
\bibinfo{author}{Schäfer, M.}, \bibinfo{author}{Strohmeier, M.},
  \bibinfo{author}{Lenders, V.}, \bibinfo{author}{Martinovic, I.},
  \bibinfo{author}{Wilhelm, M.}, \bibinfo{year}{2014}.
\newblock \bibinfo{title}{Bringing up {OpenSky}: {A} large-scale {ADS}-{B}
  sensor network for research}, \bibinfo{publisher}{IEEE}. pp.
  \bibinfo{pages}{83--94}.
\bibitem[{Shaw(1993)}]{shaw_hub_1993}
\bibinfo{author}{Shaw, S.L.}, \bibinfo{year}{1993}.
\newblock \bibinfo{title}{Hub structures of major {US} passenger airlines}.
\newblock \bibinfo{journal}{Journal of Transport Geography}
  \bibinfo{volume}{1}, \bibinfo{pages}{47--58}.
\newblock \URLprefix
  \url{https://www.sciencedirect.com/science/article/pii/096669239390037Z},
  \DOIprefix\doi{10.1016/0966-6923(93)90037-Z}.
\bibitem[{Topcuoglu et~al.(2005)Topcuoglu, Corut, Ermis and
  Yilmaz}]{topcuoglu_solving_2005}
\bibinfo{author}{Topcuoglu, H.}, \bibinfo{author}{Corut, F.},
  \bibinfo{author}{Ermis, M.}, \bibinfo{author}{Yilmaz, G.},
  \bibinfo{year}{2005}.
\newblock \bibinfo{title}{Solving the uncapacitated hub location problem using
  genetic algorithms}.
\newblock \bibinfo{journal}{Computers \& Operations Research}
  \bibinfo{volume}{32}, \bibinfo{pages}{967--984}.
\newblock \URLprefix
  \url{https://www.sciencedirect.com/science/article/pii/S030505480300279X},
  \DOIprefix\doi{10.1016/j.cor.2003.09.008}.
\bibitem[{{U.S. Congress}(2024)}]{us_congress_49_2024}
\bibinfo{author}{{U.S. Congress}}, \bibinfo{year}{2024}.
\newblock \bibinfo{title}{49 {U}.{S}. {Code} § 47102 - {Definitions}}.
\newblock \URLprefix \url{https://www.law.cornell.edu/uscode/text/49/47102}.
  \bibinfo{note}{united States Code}.
\bibitem[{{U.S. Department of
  Transportation}(2019)}]{us_department_of_transportation_alternate_2019}
\bibinfo{author}{{U.S. Department of Transportation}}, \bibinfo{year}{2019}.
\newblock \bibinfo{title}{Alternate {Essential} {Air} {Service}}.
\newblock \URLprefix
  \url{https://www.transportation.gov/office-policy/aviation-policy/alternate-essential-air-service}.
\bibitem[{{U.S. Department of
  Transportation}(2023)}]{us_department_of_transportation_dot_2023}
\bibinfo{author}{{U.S. Department of Transportation}}, \bibinfo{year}{2023}.
\newblock \bibinfo{title}{{DOT} {Penalizes} {Southwest} {Airlines} \$140
  {Million} for 2022 {Holiday} {Meltdown} {\textbar} {US} {Department} of
  {Transportation}}.
\newblock \URLprefix
  \url{https://www.transportation.gov/briefing-room/dot-penalizes-southwest-airlines-140-million-2022-holiday-meltdown}.
\bibitem[{{U.S. Department of
  Transportation}(2024)}]{us_department_of_transportation_essential_2024}
\bibinfo{author}{{U.S. Department of Transportation}}, \bibinfo{year}{2024}.
\newblock \bibinfo{title}{Essential {Air} {Service}}.
\newblock \URLprefix
  \url{https://www.transportation.gov/policy/aviation-policy/small-community-rural-air-service/essential-air-service}.
\bibitem[{Wang and Gao(2021)}]{wang_literature_2021}
\bibinfo{author}{Wang, S.}, \bibinfo{author}{Gao, Y.}, \bibinfo{year}{2021}.
\newblock \bibinfo{title}{A literature review and citation analyses of air
  travel demand studies published between 2010 and 2020}.
\newblock \bibinfo{journal}{Journal of Air Transport Management}
  \bibinfo{volume}{97}, \bibinfo{pages}{102135}.
\newblock \URLprefix
  \url{https://www.sciencedirect.com/science/article/pii/S0969699721001174},
  \DOIprefix\doi{10.1016/j.jairtraman.2021.102135}.
\bibitem[{Wei and Grubesic(2015)}]{wei_typology_2015}
\bibinfo{author}{Wei, F.}, \bibinfo{author}{Grubesic, T.H.},
  \bibinfo{year}{2015}.
\newblock \bibinfo{title}{A typology of rural airports in the {United}
  {States}: {Evaluating} network accessibility}.
\newblock \bibinfo{journal}{Review of Regional Studies} \bibinfo{volume}{45},
  \bibinfo{pages}{57--85}.
\bibitem[{Weiszer et~al.(2025)Weiszer, Delgado, Menendez-Pidal and
  Solutions}]{weiszer_air-rail_2025}
\bibinfo{author}{Weiszer, M.}, \bibinfo{author}{Delgado, L.},
  \bibinfo{author}{Menendez-Pidal, L.}, \bibinfo{author}{Solutions, N.},
  \bibinfo{year}{2025}.
\newblock \bibinfo{title}{Air-rail multimodal disruption management}.
\newblock \bibinfo{journal}{SESAR Innovation Days} .
\bibitem[{Weitz et~al.(2024)Weitz, Li, Frewin, Ellis, Sgorcea, Park, Stanley
  and Taylor}]{weitz_nas_2024}
\bibinfo{author}{Weitz, L.A.}, \bibinfo{author}{Li, M.Z.},
  \bibinfo{author}{Frewin, A.}, \bibinfo{author}{Ellis, L.},
  \bibinfo{author}{Sgorcea, R.}, \bibinfo{author}{Park, Y.},
  \bibinfo{author}{Stanley, R.}, \bibinfo{author}{Taylor, C.P.},
  \bibinfo{year}{2024}.
\newblock \bibinfo{title}{{NAS} {Traffic} {Flow} {Management}
  {Decision}-{Making} {Capability} with {Network}-{Wide} {Resilience}
  {Metrics}}, in: \bibinfo{booktitle}{{AIAA} {SCITECH} 2024 {Forum}}.
  \bibinfo{publisher}{American Institute of Aeronautics and Astronautics}.
\newblock \URLprefix \url{https://arc.aiaa.org/doi/abs/10.2514/6.2024-0573},
  \DOIprefix\doi{10.2514/6.2024-0573}. \bibinfo{note}{\_eprint:
  https://arc.aiaa.org/doi/pdf/10.2514/6.2024-0573}.
\bibitem[{Xu et~al.(2025)Xu, Hansen and Ryerson}]{xu_identification_2025}
\bibinfo{author}{Xu, J.}, \bibinfo{author}{Hansen, M.},
  \bibinfo{author}{Ryerson, M.}, \bibinfo{year}{2025}.
\newblock \bibinfo{title}{Identification and {Characterization} for
  {Disruptions} in the {U}.{S}. {National} {Airspace} {System} ({NAS})}, in:
  \bibinfo{booktitle}{First {US}-{Europe} {Air} {Transportation} {Research} \&
  {Development} {Symposium}}.
\newblock \URLprefix \url{http://arxiv.org/abs/2502.18687}.
\bibitem[{Xu et~al.(2023)Xu, Wandelt and Sun}]{xu_immuner_2023}
\bibinfo{author}{Xu, Y.}, \bibinfo{author}{Wandelt, S.}, \bibinfo{author}{Sun,
  X.}, \bibinfo{year}{2023}.
\newblock \bibinfo{title}{{IMMUNER}: {Integrated} {Multimodal} {Mobility}
  {Under} {Network} {Disruptions}}.
\newblock \bibinfo{journal}{IEEE Transactions on Intelligent Transportation
  Systems} \bibinfo{volume}{24}, \bibinfo{pages}{1480--1494}.
\newblock \URLprefix
  \url{https://ieeexplore.ieee.org/abstract/document/10027853/citations},
  \DOIprefix\doi{10.1109/TITS.2022.3224413}.
\bibitem[{Zhang et~al.(2015)Zhang, Miller-Hooks and
  Denny}]{zhang_assessing_2015}
\bibinfo{author}{Zhang, X.}, \bibinfo{author}{Miller-Hooks, E.},
  \bibinfo{author}{Denny, K.}, \bibinfo{year}{2015}.
\newblock \bibinfo{title}{Assessing the role of network topology in
  transportation network resilience}.
\newblock \bibinfo{journal}{Journal of Transport Geography}
  \bibinfo{volume}{46}, \bibinfo{pages}{35--45}.
\newblock \URLprefix
  \url{https://linkinghub.elsevier.com/retrieve/pii/S0966692315000794},
  \DOIprefix\doi{10.1016/j.jtrangeo.2015.05.006}.
\bibitem[{Zhang and Hansen(2008)}]{zhang_real-time_2008}
\bibinfo{author}{Zhang, Y.}, \bibinfo{author}{Hansen, M.},
  \bibinfo{year}{2008}.
\newblock \bibinfo{title}{Real-{Time} {Intermodal} {Substitution}: {Strategy}
  for {Airline} {Recovery} from {Schedule} {Perturbation} and for {Mitigation}
  of {Airport} {Congestion}}.
\newblock \bibinfo{journal}{Transportation Research Record}
  \bibinfo{volume}{2052}, \bibinfo{pages}{90--99}.
\newblock \URLprefix \url{https://doi.org/10.3141/2052-11},
  \DOIprefix\doi{10.3141/2052-11}.

\end{thebibliography}

\end{document}